\documentstyle[12pt]{article}
 \input amssym.def
\input amssym
\begin {document}
\def \Z{\Bbb Z}
\def \C{\Bbb C}

\def \Q{\Bbb Q}

\def \End{{\rm End}}

\def \Irr {{\rm Irr}}

\def \Hom{{\rm Hom}}

\def \<{\langle}
\def \>{\rangle}

\def \a{\alpha }

\def \L{\Lambda }

\def \b{\beta }

\def \be{\begin{equation}\label}
\def \ee{\end{equation}}
\def \bl{\begin{lem}\label}
\def \el{\end{lem}}
\def \bt{\begin{thm}\label}
\def \et{\end{thm}}
\def \bp{\begin{prop}\label}
\def \ep{\end{prop}}
\def \br{\begin{rem}\label}
\def \er{\end{rem}}
\def \bc{\begin{coro}\label}
\def \ec{\end{coro}}
\def \bd{\begin{de}\label}
\def \ed{\end{de}}
\def \pf{{\bf Proof. }}
\def \voa{{vertex operator algebra}}

\newtheorem{thm}{Theorem}[section]
\newtheorem{prop}[thm]{Proposition}
\newtheorem{coro}[thm]{Corollary}
\newtheorem{conj}[thm]{Conjecture}
\newtheorem{lem}[thm]{Lemma}
\newtheorem{rem}[thm]{Remark}
\newtheorem{de}[thm]{Definition}

\makeatletter
\@addtoreset{equation}{section}

\makeatother
\makeatletter

\baselineskip=16pt
\begin{center}{\Large \bf Simple currents and extensions of vertex operator
algebras}

\vspace{0.5cm}
Chongying Dong\footnote{Supported by NSF grant DMS-9303374 and a
research grant from the Committee on Research, UC Santa Cruz.},
Haisheng Li and Geoffrey Mason\footnote{Supported by NSF grant
DMS-9401272 and a research grant from the Committee on Research, UC
Santa Cruz.}\\ Department of Mathematics, University of California,
Santa Cruz, CA 95064
\end{center}

{\bf Abstract} We consider how a vertex operator algebra
can be extended to an abelian intertwining algebra
by a family of weak twisted modules which are {\em simple currents} associated
with semisimple weight one primary vectors. In the case that the extension
is again a vertex operator algebra, the rationality of the
extended algebra is discussed. These results are applied
to affine Kac-Moody algebras in order to
construct all the simple currents explicitly
(except for $E_8$) and to
get various extensions of the vertex operator algebras associated with
integrable representations.

\section{Introduction}

Introduced in [B] and [FLM], {\em vertex operator algebras}
are essentially {\em chiral algebras} as formulated in
[BPZ] and [MoS], and provide a powerful algebraic tool for
studying the general structure of conformal field theory. For a vertex
operator algebra $V$, one wishes to adjoin certain simple
$V$-modules to get a larger algebraic structure so that certain data such
as fusion rules and braiding matrices are naturally incorporated.
The introduction of the notions of {\em generalized vertex
(operator) algebra} and {\em abelian intertwining
algebra} in [DL] was made in this spirit. A similar notion called vertex
operator {\em para-algebra} was independently introduced and studied
in [FFR] with different motivations. Also see [M].

In this paper, we study how a vertex operator algebra
can be extended to an abelian intertwining algebra
by a family of weak twisted modules which are {\em simple currents} associated
with semisimple weight one primary vectors. In the case that the extension
is again a vertex operator algebra, we discuss the rationality of the
extended algebra. Applying these results to affine Kac-Moody algebras
we construct all the simple currents explicitly (except for $E_8$) and
get various extensions of the vertex operator algebras associated with
integrable representations. Many of the ideas discussed here are
natural continuations of ideas discussed in [DL] and [Li4]. Recently,
it is shown in [Hua] that the extension of the moonshine module \voa\
$V^{\natural}$ [FLM] by $\Z_2$-twisted module forms an abelian intertwining
algebra, where such an extension is called a nonmeromorphic extension.

Let $G$ be a torsion group of automorphisms of
$V$  and $g\in G.$
A $G$-simple current
for a vertex operator algebra is an irreducible
weak $g$-twisted module which gives a bijection between
the  equivalence
classes of irreducible weak $h$-twisted modules and the equivalence classes
of irreducible weak $gh$-twisted modules under
the ``tensor product'' if $h\in G$ commutes with $g.$ The concept
of simple current was originally introduced in [SY] for constructing
modular invariant partition functions (see also [FG] and [GW]).
A {\em simple} vertex operator algebra
$V$ is always a $G$-simple current for any $G.$
It is natural to expect that any simple
current can be {\em deformed} from $V$ by introducing a new action. Via
this principle, a class of simple currents are constructed in the present
paper. Let $U(V[1])$ be
the universal enveloping algebra of the Lie algebra
$V[1]=V\otimes\C[t,t^{-1}]/D(V\otimes\C[t,t^{-1}])$ defined in [Li3] (for
more detail see Section 2) and
let $\Delta(z)\in U(V[1])\{z\}$ satisfy conditions (\ref{li1})-(\ref{li4})
below. Then for any weak $h$-twisted $V$-module $(M,Y_{M}(\cdot,z))$, we show
that $(\tilde{M},Y_{\tilde{M}}(\cdot,z)):= (M,Y_{M}(\Delta(z)\cdot,z))$ is
a weak $gh$-twisted $V$-module where $h$ is any automorphism of $V$ of
finite order which commutes with $g$ and is not necessarily in $G.$
If $\Delta(z)\in U(V^G[1])$ ($V^G$ is the space of $G$-invariants of $V$)
is invertible and $h\in G$, then $\tilde{M}$ is isomorphic
to a tensor product module of $M$ with $\tilde{V}$ and thus
$\tilde{V}$ is a $G$-simple current. There is a simple way to construct
such $\Delta(z)$ associated to any weight one primary vector $\alpha$
of $V$ whose
component operator $\alpha(0)$ is semisimple on $V.$ The corresponding
automorphism of $V$ is given by $e^{2\pi i \alpha(0)}.$ These results
have been obtained in [Li4] in the case that $g=1.$

This paper is organized as follows: In Section 2 we recall the definitions
of weak twisted modules from [D2] and [FFR] and of intertwining operators among
weak twisted modules from [FHL] and [X] using the language of formal variables.
We
present a notion of tensor product of two weak twisted modules in terms
of universal mapping properties and relate the fusion rules (which are the
dimensions of the space of intertwining operators of certain types)
to the dimensions of certain spaces of homomorphisms of weak
twisted modules. For certain elements $\Delta(z)$ in $U(V[1])\{z\}$ associated
to
an automorphism $g$ of $V$ of finite order we show how a weak $h$-twisted
module can be deformed to a weak $gh$-twisted module if the automorphism
$h$ commutes with $g.$ We discuss the relations among this deformation,
the intertwining operators and tensor product. It turns out that
the deformation of $V$ is always a simple current. Such $\Delta(z)$
are constructed for semisimple primary vectors of $V$ of weight one.
These results
are interpreted in the theory of vertex operator algebras associated
with even lattices at the end of this section.

Section 3 deals with the extension of a simple vertex operator
algebra by a family of simple currents constructed from semisimple
primary vectors of weight one. We begin with a finite dimensional
subspace $H$ of $V_1$ which contains all semisimple vectors, and a lattice
$L$ contained in
$H$ such that the component operators $\alpha(0)$ ($\alpha\in L$)
have only rational eigenvalues on $V.$ We show that the direct sum $U$
of all the deformations of $V$ by using $\Delta(z)$ associated with
$\alpha\in L$ is a generalized vertex algebra in the sense of [DL].
Then there is an even sublattice $L_0$ such that the corresponding
deformation of $V$ is isomorphic to $V.$
We then prove that the quotient $\bar U$ of $U$ modulo the isomorphic
relations has a structure of abelian intertwining algebra. This
section is technically complicated, involving the abelian cohomology
of abelian groups introduced by Eilenberg-MacLane. We refer
the reader to [DL] for the definition of abelian intertwining algebras
and related terminology.

In Section 4 we discuss the rationality of a simple vertex operator
algebra $V=\oplus_{g\in G}V^g$ which is graded by a finite abelian group
$G$ satisfying certain conditions. Such vertex operator algebras
arise naturally in Section 3 and in [DL]. Under a mild assumption
we show how the rationality of $V^0$ implies the rationality of $V.$
In particular, if $G$ consists of automorphisms constructed from
semisimple primary weight one vectors, the assumption always hold.
In Section 5 we apply the results obtained in the previous sections
to affine algebras. After discussing simple currents
for the vertex operator algebras associated to integrable
representations and which are related to the work of [FG],
we obtain various extensions of these vertex operator algebras
by simple currents. Another result in this section is about the
representations of these algebras.
Under the assumption that certain elements in the Heisenberg subalgebra
of a given affine algebra act nilpotently on a weak module for
the vertex operator algebra we show the complete reducibility of this weak
module\footnote{This assumption has been removed recently in [DLM1].}.
In particular, this shows that any irreducible weak module is a standard
module.

\section{Simple currents and twisted modules}

In this section we first recall the definitions of twisted modules
(cf. [D2] and [FFR]) and intertwining operators among twisted
modules (cf. [FHL] and [X]). We then discuss how a weak module for a vertex
operator algebra $V$
can be deformed to a twisted module by using certain elements in the
vector space of formal power series with coefficients in the universal
enveloping algebra $U(V[1])$ of $V[1]$ which is defined below.
In general, we exhibit
the deformation of intertwining operators among weak $V$-modules to
intertwining operators among weak twisted modules.
We also apply these results to vertex operator algebras associated with
even positive-definite lattices. The ideas and
techniques used in this section derive from those in [Li4].

Let $(V,Y,{\bf 1},\omega)$ be a vertex operator
algebra (cf. [B], [FHL] and [FLM])
and let $g$ be an automorphism of $V$ of finite order $T.$  Then $V$ is a
direct
sum of eigenspaces of $g:$
$$V=\oplus_{r\in \Z/T\Z}V^r$$
where $V^r=\{v\in V|gv=e^{2\pi ir/T}v\}.$
(We abuse notation and use $ r \in \{0, 1,\ldots, T-1\}$ to denote both
an integer and the corresponding residue class.)
Following [D2] and [FFR], a weak $g$-twisted
$module$ $M$ for $V$ is a vector space equipped with a linear map
$$\begin{array}{l}
V\to (\mbox{End}\,M)\{z\}\label{map}\\
v\mapsto\displaystyle{ Y_M(v,z)=\sum_{n\in\Q}v_nz^{-n-1}\ \ \ (v_n\in
\mbox{End}\,M)}
\end{array}$$
(where for any vector space $W,$ we define $W\{z\}$ to be the vector
space of $W$-valued formal series in $z,$ with arbitrary complex powers of $z$)
satisfying the following conditions for $u,v\in V$,
$w\in M$, and $r\in\Z/T\Z$:
\begin{eqnarray}
& &Y_M(v,z)=\sum_{n\in \frac{r}{T}+\Z}v_nz^{-n-1}\ \ \ \ \mbox{for}\ \
v\in V^r;\label{1/2}\\
& &v_lw=0\ \ \
\mbox{for}\ \ \ l\in \Q \ \ \mbox{sufficiently\ large};\label{vlw0}\\
& &Y_M({\bf 1},z)=1;\label{vacuum}
\end{eqnarray}
 \begin{equation}\label{jacobi}
\begin{array}{c}
\displaystyle{z^{-1}_0\delta\left(\frac{z_1-z_2}{z_0}\right)
Y_M(u,z_1)Y_M(v,z_2)-z^{-1}_0\delta\left(\frac{z_2-z_1}{-z_0}\right)
Y_M(v,z_2)Y_M(u,z_1)}\\
\displaystyle{=z_2^{-1}\left(\frac{z_1-z_0}{z_2}\right)^{-r/T}
\delta\left(\frac{z_1-z_0}{z_2}\right)
Y_M(Y(u,z_0)v,z_2)}
\end{array}
\end{equation}
if $u\in V^r;$
$$[L(m),L(n)]=(m-n)L(m+n)+\frac{1}{12}(m^3-m)\delta_{m+n,0}(\mbox{rank}\,V)$$
for $m, n\in {\Z},$ where
\begin{eqnarray}
& &L(n)=\omega_{n+1}\ \ \ \mbox{for}\ \ \ n\in{\Z}, \ \ \
\mbox{i.e.},\ \ \ Y_M(\omega,z)=\sum_{n\in{\Z}}L(n)z^{-n-2};\nonumber\\
& &\frac{d}{dz}Y_M(v,z)=Y_M(L(-1)v,z).\label{6.72}
\end{eqnarray}
This completes the definition. We denote this module by
$(M,Y_M)$ (or briefly by $M$).

\begin{rem}\label{r2.0} If a weak $g$-twisted $V$-module $M$ further
is a $\C$-graded
vector space:
$$M=\coprod_{\lambda \in{\C}}M_{\lambda},$$
such that for each $\lambda \in \C $
$$\dim M_{\lambda}<\infty\;$$
$$M_{\lambda+\frac{n}{T}}=0$$
for $n\in {\Bbb Z}$ sufficiently small, and
$$L(0)w=nw=(\mbox{wt}\,w)w \ \ \ \mbox{for}\ \ \ w\in M_n\
(n\in{\C}),$$
we call $M$ a $g$-twisted $V$-module.
\end{rem}

A weak $1$-twisted $V$-module $M$
is called a {\em weak $V$-module} and a $1$-twisted $V$-module is called
a $V$-{\em module}. A $g$-{\em homomorphism} $f$ from a weak $g$-twisted
$V$-module $M$ to another weak $g$-twisted
$V$-module $W$ is a linear map $f: M\to W$ such that
$$fY_M(u,z)=Y_W(u,z)f$$
for all $u\in V.$ We denote the space of all $g$-homomorphisms from
$M$ to $W$ by $\Hom_g(M,W).$
A $g$-{\em isomorphism} is a bijective $g$-homorphism.

Next we shall define intertwining operators among weak $g_k$-twisted
modules $(M_k,Y_{M_k})$ for $k=1,2,3$ where $g_k$ are commuting
automorphisms of order $T_k$ (cf. [X]). In this case
$V$ decomposes into the direct sum of common
eigenspaces
$$V=\oplus_{j_1,j_2}V^{(j_1,j_2)}$$
where
$$V^{(j_1,j_2)}=\{v\in V|g_kv=e^{2\pi ij_k/T_k}, k=1,2\}$$

An $intertwining$ $operator$ $of$ $type$ {
$\left[\hspace{-3 pt}\begin{array}{c}M_3\\ M_1\
M_2\end{array}\hspace{-3 pt}\right]$}
associated with the given data is a linear map
\begin{equation}\label{11.operator}
\begin{array}{lll}
M_1&\to& (\mbox{Hom}(M_2,M_3))\{z\} \\
w&\mapsto& {\cal Y}(w,z)=\displaystyle{\sum_{n\in{\C}}w_nz^{-n-1}}
\end{array}
\end{equation}
such that for $w^i\in M_i$ ($i=1,2$),
fixed $c\in \C$ and $n\in \Q$ sufficiently large
\begin{equation}\label{a1}
w^1_{c+n}w^2=0;
\end{equation}
the following (generalized) Jacobi identity holds on $M_2:$ for
$u\in V^{(j_1,j_2)}$ and $w\in M_1,$
$$z^{-1}_0\left(\frac{z_1-z_2}{z_0}\right)^{j_1/T_1}
\delta\left(\frac{z_1-z_2}{z_0}\right)Y_{M_3}(u,z_1){\cal
Y}(w,z_2)$$
\begin{equation}\label{11.rjac}
\ \ \ \ -z_0^{-1}\left(\!\frac{z_2\!-\!z_1}{z_0}\!\right)^{j_1/T_1}\!
\delta\!\left(\!\frac{z_2\!-\!z_1}{-z_0}\!\right)\!
{\cal Y}(w,z_2)Y_{M_2}(u,z_1)
\end{equation}
$$=z_2^{-1}\left(\frac{z_1-z_0}{z_2}\right)^{-j_2/T_2}\delta
\left(\frac{z_1-z_0}{z_2}\right)
{\cal Y}(Y_{M_1}(u,z_0)w,z_2)$$
and
\begin{equation}\label{11.l-1}
\frac{d}{dz}{\cal Y}(w,z)={\cal Y}(L(-1)w,z),
\end{equation}
where $L(-1)$ is the operator acting on $M_1.$

The intertwining operators of type {
$\left[\hspace{-3
pt}\begin{array}{c}M_3\\ M_1\ M_2\end{array}\hspace{-3
pt}\right]${} }
associated with prescribed data clearly form a vector space,
which we denote by ${\cal V}^{M_3}_{M_1M_2}.$  We set
\medskip
\begin{equation}
N^{M_3}_{M_1M_2}=\mbox{dim}\,{\cal V}^{M_3}_{M_1M_2}.
\end{equation}
These numbers are called the $fusion$ $rules$ associated with the
algebra, modules and auxiliary data. It is easy to observe that
if $N^{M_3}_{M_1M_2}>0$ then $g_3=g_1g_2$ (see [X]). Thus we shall
assume this relation in the following discussion.

The fusion rules have certain symmetry properties. Our next goal is
to show that $N^{M_3}_{M_2M_1}=N^{M_3}_{M_1M_2}.$ Let ${\cal Y}$ be an
intertwining operator of type {
$\left[\hspace{-3 pt}\begin{array}{c}M_3\\ M_1\
M_2\end{array}\hspace{-3 pt}\right]$}. We define a linear map
\begin{equation}\label{11s}
\begin{array}{lll}
M_2&\to& (\mbox{Hom}(M_1,M_3))\{z\} \\
w&\mapsto& {\cal Y}^{\pm}(w,z)=\displaystyle{\sum_{n\in{\C}}w_nz^{-n-1}}
\end{array}
\end{equation}
by the skew-symmetry
\begin{equation}\label{s2}
{\cal Y}^{\pm}(w^2,z)w^1=e^{zL(-1)}{\cal Y}(w^1,e^{\pm \pi i}z)w^2
\end{equation}
for $w^i\in M_i,$ where $Y(w^1,e^{\pm\pi i}z)=\sum_{n\in
\C}w^1_ne^{\pm(-n-1)\pi i}z^{-n-1}.$
\begin{lem}\label{la1} (1)
The operator ${\cal Y}^{\pm}(\cdot,z)$ is an
intertwining operator of type {
$\left[\hspace{-3 pt}\begin{array}{c}M_3\\ M_2\
M_1\end{array}\hspace{-3 pt}\right]$}.

(2) Each of the two maps ${\cal Y}\mapsto {\cal Y}^{\pm}$ is a linear
isomorphism
{}from  ${\cal V}^{M_3}_{M_1M_2}$ to ${\cal V}^{M_3}_{M_2M_1}.$ In particular
$N^{M_3}_{M_2M_1}=N^{M_3}_{M_1M_2}.$
\end{lem}

\pf (1) The relation (\ref{a1}) is clear. In order to prove (\ref{11.l-1})
for ${\cal Y}^{\pm}(w^2,z)$ we first
observe from (\ref{11.rjac}) that
$$[L(-1),{\cal Y}(w^1,z)]={\cal Y}(L(-1)w^1,z)=\frac{d}{dz}{\cal Y}(w^1,z),$$
$$e^{L(-1)z_0}{\cal Y}(L(-1)w^1,z)e^{-L(-1)z_0}={\cal Y}(L(-1)w^1,z+z_0).$$
Thus we have the following calculation:
\begin{eqnarray*}
& &\frac{d}{dz}{\cal Y}^{\pm}(w^2,z)w^1=\frac{d}{dz}e^{zL(-1)}{\cal
Y}(w^1,e^{\pm i\pi}z)w^2\\
& &\ \ \ =e^{zL(-1)}L(-1){\cal Y}(w^1,e^{\pm i\pi}z)w^2
          -e^{zL(-1)}{\cal
Y}(L(-1)w^1,e^{\pm i\pi}z)w^2\\
& &\ \ \ =e^{zL(-1)}{\cal Y}(w^1,e^{\pm i\pi}z)L(-1)w^2\\
& &\ \ \ ={\cal Y}^{\pm}(L(-1)w^2,z)w^1.
\end{eqnarray*}
The proof of Proposition 2.2.2 of [G] provides (as a special case)
a proof of the Jacobi identity  (\ref{11.rjac}) for $Y_{M_i}(u,z_1)$ and
${\cal Y}^{\pm}(w^2,z_2).$

(2) An easy verification shows that ${\cal Y}^{+-}={\cal Y}$ for any
${\cal Y}\in {\cal V}^{M_3}_{M_1M_2}.$ The isomorphism between
${\cal V}^{M_3}_{M_1M_2}$ and ${\cal V}^{M_3}_{M_2M_1}$ is now clear.
\ \ $\Box$
\bigskip

Next we formulate the notion of tensor product
$M_1 \boxtimes M_2$ of $M_1$ and $M_2$: it is a weak $g_3$-twisted $V$-module
defined by a universal mapping property. See [Li3] and [HL1-HL2] for
the definition of tensor product of (ordinary) weak modules.

\begin{de}\label{d1}A
 {\it tensor product} for the ordered pair $(M_1,M_{2})$ is a pair
$(M,F(\cdot,z))$ consisting of a weak $g_3$-twisted
$V$-module $M$ and an intertwining
operator $F(\cdot,z)$ of type
$\left(\begin{array}{c}M\\M_{1}\,M_{2}\end{array}\right)$ such that
the following universal
property holds: for any weak $g_3$-twisted $V$-module $W$ and any intertwining
operator $I(\cdot,z)$ of type
$\left(\begin{array}{c}W\\M_{1}\,M_2\end{array}\right)$, there
exists a unique $V$-homomorphism $\psi$ from $M$ to $W$ such that
$I(\cdot,z)=\psi\circ F(\cdot,z)$. (Here $\psi$ extends canonically to a
linear map from $M\{z\}$ to $W\{z\}$.)
\end{de}

The following proposition is a direct consequence of the definition and
Lemma \ref{la1}.
\begin{prop}\label{p2.7} Let $(M,F(\cdot,z))$ be a tensor
product of $M_1$ and $M_2.$

(1) For any weak $g_3$-twisted $V$-module
$M_{3}$, the space ${\rm Hom}_{g_3}(M,M_{3})$
is linearly isomorphic to ${\cal V}_{M_1M_2}^{M_3}.$

(2) The pair $(M,F^{\pm}(\cdot,z))$ is a tensor product of $M_2$ and
$M_1.$ In particular, the tensor product is commutative.

(3) If $g_1=1$ and $M_1=V$ then $M$ is isomorphic to $M_2.$ That is,
$V \boxtimes M_2=M_2.$
\end{prop}

Let $V$ be a vertex operator algebra and $g$ an automorphism of $V$ of order
$T.$ We recall from [DLM2] the Lie algebra
$$V[g]=\oplus_{i=0}^{T-1}(V^i\otimes
t^{i/T}\C[t,t^{-1}]/D(\oplus_{i=0}^{T-1}(V^i\otimes t^{i/T}\C[t,t^{-1}]))$$
 associated with $V$ and $g,$ with bracket
$$[u(m),v(n)]=\sum_{i=0}^{\infty}{m\choose i}(u_iv)(m+n-i-1),$$
where $D=L(-1)\otimes 1+\frac{d}{dt}\otimes 1$ and $u(m)$ is the image of
$u\otimes t^m$ in $V[g].$ Then $V^0[1]$ is a Lie subalgebra of both
$V[1]$ and $V[g],$ and $V[g]$ acts on any weak $g$-twisted $V$-module.

Let $\Delta(z)\in U(V[1])\{z\}$ satisfy the following conditions:
\begin{eqnarray}
& &z^{{j\over T}}\Delta(z)a\in V[z,z^{-1}] \mbox{{\it for }}a\in V^j;
\label{li1}\\
& &\Delta(z){\bf 1}={\bf 1};\label{li2}\\
& &[L(-1), \Delta(z)]=- {d\over dz}\Delta(z);\label{li3}\\
& &Y(\Delta(z_{2}+z_{0})a,z_{0})\Delta(z_{2})
=\Delta(z_{2})Y(a,z_{0})\;\;\;\mbox{{\it for any }}a\in V.\label{li4}
\end{eqnarray}
Let $G(V,g)$ be the
set of all $\Delta(z)$ satisfying the conditions (\ref{li1})-(\ref{li4}).
Define $G^0(V,g)$ to be these $\Delta(z)$ in $G(V,g)$ which are invertible.
The following Lemma is obvious:
\begin{lem}\label{l2.0} (1) {\it Let $\Delta_{1}(z), \Delta_{2}(z)\in G(V,g)$.
Then
$\Delta_{1}(z)\Delta_{2}(z)\in G(V,g^{2})$.} In particular.
$G(V,1)$ is a semigroup with $id_{V}\in G(V,1).$

(2) {\it If $\Delta(z)\in G(V,g)$ has an inverse $\Delta^{-1}(z)\in
U(V[1])\{z\}$. Then $\Delta^{-1}(z)\in G(V, g^{-1})$.}
\end{lem}

{}From now on we fix two commutative automorphisms $g$ and $h$ of $V$ of order
$S$ and $T$ respectively.
The following result generalizes the corresponding results
in [Li4] with $g=h=1.$
\begin{lem}\label{l2.1}
Let $(M,Y_{M}(\cdot,z))$ be a weak $h$-twisted $V$-module and $\Delta(z)\in
G(V,g).$
Set $\tilde{M}=M$ and $Y_{\tilde{M}}(\cdot,z)=Y_{M}(\Delta(z)\cdot,z)$. Then
$(\tilde{M},Y_{\tilde{M}}(\cdot,z))$ is a weak $gh$-twisted
$V$-module.
\end{lem}

{\bf Proof.} First, (\ref{li2}) implies that $\tilde{Y}_{\tilde{M}}({\bf
1},z)=Id_{\tilde{M}}$.
Second, it is easy to see that (\ref{li3}) implies that
$\tilde{Y}_{\tilde{M}}(L(-1)a,z)={d\over dz}\tilde{Y}_{\tilde{M}}(a,z)$ for any
$a\in V$.
Let $a\in V^{(i,j)}$ and $b\in V$. Then we have:
\begin{eqnarray}
& &\ \ z_{0}^{-1}\delta\left(\frac{z_{1}-z_{2}}{z_{0}}\right)
Y_{M}(\Delta(z_{1})a,z_{1})Y_{M}(\Delta(z_{2})b,z_{2})\nonumber\\
& &\ \ -z_{0}^{-1}\delta\left(\frac{z_{2}-z_{1}}{-z_{0}}\right)
Y_{M}(\Delta(z_{2})b,z_{2})Y_{M}(\Delta(z_{1})a,z_{1})\nonumber\\
& &=z_{2}^{-1}\left(\frac{z_{1}-z_{0}}{z_{2}}\right)^{\frac{-i}{S}}
\delta\left(\frac{z_{1}-z_{0}}{z_{2}}\right)
Y_{M}(Y(\Delta(z_{1})a,z_{0})\Delta(z_{2})b,z_{2})\nonumber\\
& &=z_{2}^{-1}\left(\frac{z_{1}-z_{0}}{z_{2}}\right)^{\frac{-i}{S}}
\delta\left(\frac{z_{1}-z_{0}}{z_{2}}\right)z_{1}^{-{j\over T}}
Y_{M}(Y((z_{2}+z_{0})^{{j\over
T}}\Delta(z_{2}+z_{0})a,z_{0})\Delta(z_{2})b,z_{2})\nonumber\\
& &=z_{2}^{-1}\left(\frac{z_{1}-z_{0}}{z_{2}}\right)^{\frac{-i}{S}}
\delta\left(\frac{z_{1}-z_{0}}{z_{2}}\right)
\left(\frac{z_{2}+z_{0}}{z_{1}}\right)^{{j\over T}}
Y_{M}(Y(\Delta(z_{2}+z_{0})a,z_{0})\Delta(z_{2})b,z_{2})\nonumber\\
& &=z_{2}^{-1}\left(\frac{z_{1}-z_{0}}{z_{2}}\right)^{\frac{-i}{S}-{j\over
T}}\delta\left(\frac{z_{1}-z_{0}}{z_{2}}\right)
Y_{M}(\Delta(z_{2})Y(a,z_{0})b,z_{2}).
\end{eqnarray}
Thus
$(\tilde{M},Y_{\tilde{M}}(\cdot,z))$ is a weak $gh$-twisted
$V$-module. \ \ $\Box$


We shall use the notation $V^h$ for the $h$-fixed-point vertex operator
subalgebra of $V$ for an automorphism $h.$
\begin{lem}\label{l2.2} (1) Assume that
$g$ commutes with each $g_i.$ {\it Let $\Delta(z)\in G(V^{g_1},g)$, let $M_{i}$
$(i=1,2,3)$ be weak $g_i$-twisted $V$-modules and
$I(\cdot,z)$ an intertwining operator of type
$\left(\begin{array}{c}M_{3}\\M_{1} M_{2}\end{array}\right)$. Then
$\tilde{I}(\cdot,z)=I(\Delta(z)\cdot,z)$ is an intertwining operator of type
$\left(\begin{array}{c}\tilde{M}_{3}\\M_{1} \tilde{M}_{2}
\end{array}\right)$.}

(2) {\it Let $\Delta(z)\in G(V,g)$ and let $\psi$ be a
$h$-homomorphism from $W$ to $M$. Then
$\psi$ is a $gh$-homomorphism from $\tilde{W}$ to
$\tilde{M}$.}

(3) Let
$\Delta(z)\in G(V^g,1)$ be such that $(V,Y(\Delta(z)\cdot,z)$ is isomorphic to
the
adjoint module $(V,Y(\cdot,z))$. Then there is
a nonzero homomorphism of weak $g$-twisted modules
{}from $(M,Y_{M}(\cdot,z))$ to $(M,Y_{M}(\Delta(z)\cdot,z))$
for any weak
$g$-twisted $V$-module $(M,Y_{M}(\cdot,z))$. Moreover, if
$\Delta(z)\in G^0(V^g,1),$ any such homomorphism is an isomorphism.
\end{lem}

\pf The proof of (1) is similar to that of Lemma \ref{l2.1} and we omit
details. (2) is a special case of (1) with $M_1=V, M_2=W$ and $M_3=M.$
It remains to show (3).
Clearly $Y_M(\cdot,z)$ is an intertwining operator of type
$\left(\begin{array}{c}M\\V
M\end{array}\right)$. By Lemma \ref{la1}, there is a nonzero intertwining
operator of type $\left(\begin{array}{c}M\\M
V\end{array}\right).$ Now by (1) there is a nonzero intertwining
operator of type $\left(\begin{array}{c}\tilde{M}\\M
\tilde{V}\end{array}\right),$ which yields a  nonzero intertwining
operator of type $\left(\begin{array}{c}\tilde{M}\\M
V\end{array}\right)$ by hypothesis. Consequently there is
a nonzero intertwining
operator $I(\cdot,z)$
of type $\left(\begin{array}{c}\tilde{M}\\V M\end{array}\right)$ by
Lemma \ref{la1} once more. Now $I(1,z)$ is the
desired  nonzero homomorphism. The other
assertions are clear. \ \ $\Box$

\begin{prop}\label{p2.8} {\it Let $(W,F(\cdot,z))$ be a tensor product
of a weak $g_1$-twisted module $M_{1}$ and a weak $g_2$-twisted
module $M_{2}$ and assume that $g$ commutes with each $g_i.$
Then if $\Delta(z)\in G^{0}(V^{g_1},g),$
$(\tilde{W},\tilde{F}(\cdot,z))$ is a tensor product of the pair
$(M_{1}, \tilde{M}_{2})$.}
\end{prop}

{\bf Proof.} First by Lemma \ref{l2.2} (1),
we have an intertwining operator
$\tilde{F}(\cdot,z)\!=\!F(\Delta(z)\cdot,z)$ of type
$\left(\begin{array}{c}\tilde{W}\\M_{1} \tilde{M}_{2}\end{array}\right)$.
Let $M$ be a weak $gg_1g_2$-twisted $V$-module and let $I(\cdot,z)$ be
any intertwining operator of type
$\left(\!\begin{array}{c}M\\M_{1} \tilde{M}_{2}\end{array}\!\right)$. Then
$I(\Delta(z)^{-1}\cdot,z)$ is an intertwining operator of type
$\left(\!\begin{array}{c}\hat{M}\\M_{1} M_{2}\end{array}\!\right)$,
where
$(\hat{M},Y_{\hat{M}}(\cdot,z))=(M,Y_{M}(\Delta(z)^{-1}\cdot,z))$. By the
universal property of $(W,
F(\cdot,z))$, there is a unique $g_1g_2$-homomorphism $\psi$ from
$W$ to $\hat{M}$
such that $\hat{I}(\cdot,z)=\psi\circ \hat{F}(\cdot,z)$. By Lemma \ref{l2.2}
(2),
$\psi$ is a $g_1g_2$-homomorphism from $W$ to $M$. Since $\Delta(z)u$ only
involves finitely many terms, we have: $I(\cdot,z)=\psi\circ
F(\cdot,z)$. It is easy to check the uniqueness, so that the proof is
complete.$\;\;\;\;\Box$

\begin{coro}\label{c2.9}
{\it Let $M$ be a weak $h$-twisted $V$-module and let
$\Delta(z)\in G^0(V^h,g)$. Then $\tilde{M}$ is isomorphic to the tensor
product module of $M$ with $\tilde{V}$.}
\end{coro}

{\bf Proof.} It is easy to observe that
$(M,F(\cdot,z))$ is a tensor
product of $M$ and $V$ where $F(\cdot,z)$ is the transpose
intertwining operator of $Y_{M}(\cdot,z)$ (also see Proposition \ref{p2.7}
(3)). By Proposition \ref{p2.8},
$(\tilde{M},\tilde{F}(\cdot,z))$ is a tensor product of $M$ and $\tilde{V}$.
$\;\;\;\;\Box$

\begin{rem}\label{phy} In the situation of Corollary \ref{c2.9} with
$g=h=1,$ $\tilde M$ is {\em defined} as the tensor product of
$M$ with $\tilde V$ in the physics literature (cf. [MaS]). Corollary
\ref{c2.9} asserts that this formulation coincides with our
axiomatic notion of tensor product in this special case.
\end{rem}

The following definition is essentially due to Schellekens and
Yankielowicz [SY].
\begin{de}\label{d2.10}
Let $V$ be a vertex operator algebra and  $G$ a torsion
group of automorphisms of $V.$
We denote the set of equivalence classes
 of irreducible weak $h$-twisted modules by $\Irr_h(V)$ for $h\in G.$
For convenience
we write $\Irr(V)\!=\Irr_1(V).$ An
irreducible weak $g$-twisted $V$-module $M$ for $g\in G$
is called a $G$-{\em simple current}  if the tensor functor
``$M\boxtimes \cdot$'' is
a bijection from $\Irr_h(V)$ to $\Irr_{gh}(V)$
for any $h\in G$ which commutes with $g.$ A $1$-simple current ($G=1$)
is called a simple current.
\end{de}

Clearly, a $G$-simple current $M\in \Irr(V)$ acts on $\Irr_h(V)$ for $h\in G$
as a permutation via the tensor product $M\boxtimes\cdot$ for any $h.$
\begin{prop}\label{p2.11}
{\it For any $\Delta(z)\in G^{0}(V^G,g)$, $(V,
Y(\Delta(z)\cdot,z))$ is a $G$-simple current if $V$ is a simple vertex
operator
algebra.}
\end{prop}

{\bf Proof.} By  Corollary \ref{c2.9} for any weak $h$-twisted module $M,$
$\tilde{M}=(M,Y_M(\Delta(z)\cdot,z)$ is isomorphic to the
tensor product of $M$ with $(V,Y(\Delta(z)\cdot,z)),$
and $M$ is isomorphic to the tensor product of $\tilde{M}$ with
$(V,Y(\Delta(z)^{-1}\cdot,z)).$ Thus if $M$ is irreducible so is $\tilde{M},$
$\tilde{V}\mapsto \tilde{V}\boxtimes M$ being a bijection from $\Irr_h(V)$ to
$\Irr_{gh}(V).$\ \ $\Box$

\begin{conj}\label{con2.14}{Let $V$ be a vertex operator algebra and
$G$ a group of automorphisms of finite order of $V.$ Define
$\bar{V}$ to be the direct sum of all $G$-simple currents of $V.$
Then $\bar{V}$ is an abelian intertwining algebra in the precise sense of
[DL].}
\end{conj}

In the next section  we will prove this conjecture in some special cases.

Let $V$ be a vertex operator algebra and let $\alpha\in V$
satisfying the following conditions:
\begin{eqnarray}\label{h2.11}
L(n)\alpha=\delta_{n,0}\alpha,\;
\alpha(n)\alpha=\delta_{n,1}\gamma {\bf 1}\;\mbox{for any }n\in {\Z}_{+},
\end{eqnarray}
where $\gamma$ is a fixed complex number. Notice that condition (\ref{h2.11})
implies
that $\a$ is a primary vector of weight one.
Furthermore, we assume that $\a(0)$ acts semisimply on $V$ with
rational eigenvalues.
It is clear that $e^{2\pi i\a(0)}$ is an automorphism of $V$. If the
denominators of all eigenvalues of $\a(0)$ are bounded, then
 $e^{2\pi i\a(0)}$ is of finite order.
Define
\begin{eqnarray}
\Delta(\a,z)=z^{\a(0)}\exp\left(\sum_{k=1}^{\infty}\frac{\a(k)}{-k}
(-z)^{-k}\right)\in U(V^{e^{2\pi i\a(0)}}[1])\{z\}.
\end{eqnarray}
Recall the following proposition from [Li2].

\begin{prop}\label{p2.14} Let $\tau$ be a finite-order automorphism of
$V$ such that $\tau\a=\a$
 and let $(M,Y_{M}(\cdot,z))$ be any
$\tau$-twisted $V$-module. Then $(M, Y(\Delta(\alpha,z)\cdot,z))$ is
a $\sigma_{\a}\tau $-twisted weak $V$-module, where
$\sigma_{\a}=e^{-2\pi i\a(0)}$.
\end{prop}

The main part of the proof of Proposition \ref{p2.14} in [Li2]
is to establish that $\Delta(\alpha,z)$ satisfies the condition (\ref{li4}).

We end this section by discussing an example, namely vertex operator algebras
associated with even positive-definite lattices.
Let $L$ be a positive definite rational lattice with form $\langle,\rangle$
 let $L_{0}$ be an even
sublattice of $L$ and let $V=V_{L_{0}}$ be the vertex operator algebra
constructed in [B] and [FLM]. For any
$\beta\in L$, $V_{\beta+L_{0}}$ is a twisted $V_{L_{0}}$-module for the inner
automorphism $\sigma_{\beta}=e^{-2\beta(0)\pi i}$ (see [DM] and [Le]).
It is easy to see that
$V_{\beta+L_{0}}$ is isomorphic to the adjoint module $V_{L_{0}}$ if and only
if $\beta\in L_{0}$.
\begin{prop}\label{p2.15} Let $\beta\in L.$ Then
as a $\sigma_{\beta}$-twisted $V_{L_{0}}$-module,
$(V_{L_{0}},Y(\Delta(\beta,z)\cdot,z))$ is isomorphic to $V_{L_{0}+\beta}$.
\end{prop}

{\bf Proof.} First, since $\Delta(\beta,z)$ is invertible and $V_{L_{0}}$ is
a simple vertex operator algebra, it follows from Proposition \ref{p2.11}
that $(V_{L_{0}},Y(\Delta(\beta,z)\cdot,z))$ is an irreducible
weak $\sigma_{\beta}$-twisted $V_{L_{0}}$-module. For any $\a\in
H={\C}\otimes_{{\Z}}L_{0}$,
we have:
\begin{eqnarray}
\Delta(\beta,z)\a=\Delta(\beta,z)\a(-1){\bf 1}=\a+z^{-1}
\langle \beta,\a\rangle.
\end{eqnarray}
Let $\psi$ be the algebra automorphism of $U(V_{L_{0}}[1])$ such that
$\psi(Y(a,z))=Y(\Delta(\beta,z)a,z)$ for $a\in V_{L_{0}}$. Then we have:
\begin{eqnarray}
\psi(\a(n))=\a(n)+\delta_{n,0}\langle\beta,\alpha\rangle\;\;\;\mbox{for }\a\in
H, n\in {\Z}.
\end{eqnarray}
Then the action of $\psi(H(0))$ on $V_{L_{0}}$
is also semisimple with $L_{0}+\beta$ as the set of $H$-weights and
$V_{L_{0}}$ is still a completely
reducible module for the Heisenberg algebra $\psi(\tilde{H})$.
It follows from  the classification result [D1] that
$(V_{L_{0}},Y(\Delta(\beta,z)\cdot,z))$ is isomorphic to
$V_{L_{0}+\beta}$. $\;\;\;\;\Box$

Let $P$ be the dual lattice of $L.$ Then
$V_{\beta+L_{0}}$ is a $V_{L_0}$-module if $\beta\in P.$
It is proved in [D1] that
there is a 1-1
correspondence between the equivalence classes of irreducible
modules for $V_{L_{0}}$ and the  cosets of $P/L_{0}.$
More specifically, $V_{P}$ is the direct sum of all inequivalent
irreducible $V_{L_0}$-modules:
\begin{eqnarray}
V_{P}= V_{L_{0}+\beta_{1}}\oplus \cdots \oplus V_{L_{0}+\beta_{k}}
\end{eqnarray}
where $k=|P/L_{0}|$. For the vertex operator algebra $V_{L_{0}}$,
intertwining operators are explicitly constructed and fusion rules are
calculated in [DL]. Moreover, it is established
in [DL] that $V_P$ is an abelian
intertwining algebra.

For any fixed $\beta\in L$, we consider all the irreducible
$\sigma_{\beta}$-twisted
$V_{L_{0}}$-modules. It was essentially proved [D1] that any irreducible
$\sigma_{\beta}$-twisted $V_{L_{0}}$-module is isomorphic to
$V_{\beta+\beta_{i}+L_{0}}$
for some $1\le i\le k$.
It follows from Proposition \ref{p2.15} that all irreducible
$\sigma_{\beta}$-twisted
$V_{L_0}$-modules
can be obtained as $(V_{L_{0}},Y(\Delta(\gamma,z)\cdot,z))$ for $\gamma\in L$,
so that
by Proposition 2.12 all irreducible $\sigma_{\beta}$-twisted
$V_{L_{0}}$-modules are
simple currents.

\section{Abelian intertwining algebras}
In this section we consider extensions of a vertex operator algebra $V,$ whose
weight one subspaces is nonzero, by incorporating certain twisted modules.
The corresponding twist elements are automorphisms
$e^{2\pi ih(0)}$ where $h\in V_1$ is a primary
vector and $h(0)$ acts semisimply on $V.$ (There should be no confusion
between $h$ in Section 2 -- an automorphism of $V,$ and $h$ in this section --
an element in $V_1.$) We prove
that such an extension is a generalized vertex operator algebra and that
the quotient space modulo isomorphic relations is
an abelian intertwining algebra in the sense of Dong
and Lepowsky [DL]. We refer the reader to [DL] for the definitions
of generalized vertex operator algebras and of abelian intertwining algebras.
In this section we assume that $V$ is a simple \voa.

Let $h\in V$ satisfy condition (\ref{h2.11}), that is, $L(n)h=\delta_{n,0}h$
and $h(n)h=\delta_{n,1}\gamma {\bf 1}$ for some fixed complex number $\gamma.$
For any $s\in
{\Q}$, set
\begin{eqnarray}
E^{\pm}(sh,z)=\exp\left(\sum_{k=1}^{\infty}\frac{sh(\pm k)}{
k}z^{\mp k}\right).
\end{eqnarray}
Then we have:
\begin{eqnarray}
E^{+}(sh,z_{1})E^{-}(th,z_{2})=\left(1-{z_{2}\over
z_{1}}\right)^{-\gamma st} E^{-}(th,z_{2})E^{+}(sh,z_{1})
\end{eqnarray}
for $s,t\in\Q$ (cf. formula (4.3.1) of [FLM]).

Let us recall some elementary results from [Li4].
\begin{lem}\label{l3.1} Let $h\in V_1$ such that (\ref{h2.11}) holds. Then
for any $a\in V$,
\begin{eqnarray}
& &e^{z(L(1)-h(1))}e^{-zL(1)}=\exp\left(\sum_{k=1}^{\infty}
\frac{h(k)}{k}(-z)^{k}\right),\\
& &e^{z(L(-1)+h(-1))}e^{-zL(-1)}=\exp\left(\sum_{k=1}^{\infty}
\frac{h(-k)}{k}z^{k}\right),
\end{eqnarray}
\begin{eqnarray}
& &Y(E^{-}(h,z_{1})a,z_{2})=E^{-}(h,z_{1}+z_{2})E^{-}(-h,z_{2})\cdot\nonumber\\
& &\ \ \ \cdot Y(a,z_{2})z_{2}^{-h(0)}E^{+}(h,z_{2})
(z_{2}+z_{1})^{h(0)}E^{+}(-h,z_{2}+z_{1}),\label{3.5}
\end{eqnarray}
\begin{eqnarray}
E^{-}(h,z_{1})Y(a,z_{2})E^{-}(-h,z_{1})
=Y(\Delta(-h,z_{2}-z_{1})\Delta(h,z_{2})a,z_{2}).\ \ \Box\label{3.6}
\end{eqnarray}
\end{lem}

Let $H$ be a finite-dimensional subspace of $V_1$
satisfying the following conditions:
\begin{eqnarray}
L(n)h=\delta_{n,0}h,\;h(n)h'=\langle h,h'\rangle \delta_{n,1}{\bf 1}
\;\;\;\mbox{for }n\in {\Z}_{+}, h,h'\in H,\label{3.7}
\end{eqnarray}
where $\langle\cdot,\cdot\rangle$ is assumed to be a nondegenerate symmetric
bilinear form
on $H$. Then we may identify $H$ with its dual $H^{*}$.
We also assume that for any $h\in H$, $h(0)$ acts semisimply on $V$. Then
\begin{eqnarray}
V=\oplus_{\alpha\in H}V^{(0,\alpha)},\;\mbox{ where }
V^{(0,\alpha)}=\{u\in V| h(0)u=\langle \alpha,h\rangle u
\;\;\mbox{for }h\in H\}.
\end{eqnarray}
Let $L$ be a lattice in $H$ such that for each $\alpha\in L$,
$\alpha(0)$ has rational eigenvalues on $V$.
{\em From now on, we assume that there is a positive integer $T$
such that the eigenvalues of $T\alpha(0)$ on $V$ are integers.} One can show
that this holds if $V$ is finitely generated. Let $G$ be the group of
automorphisms of $V$ generated by $e^{2\pi i\alpha(0)}$ for $\alpha\in L.$
Then $G$ is an abelian torsion group. Note that
$\Delta(\alpha,z)\in G^0(V^G,\sigma_{\a})$ and $(V,Y(\Delta(\a,z)\cdot,z))$ is
a $G$-simple current by Proposition \ref{p2.11} where $\sigma_{\alpha}=e^{-2\pi
i\alpha(0)}$.

For any $\alpha\in L$ and for any weak $V$-module $(M,Y_{M}(\cdot,z))$, we have
a weak $\sigma_{\alpha}$-twisted module
$$(M^{(\alpha)}, Y_{\alpha}(\cdot,z))=(M,Y_{M}(\Delta(\alpha,z)\cdot,z)).$$
This yields
a linear isomorphism $\psi_{\alpha}$ from
$M^{(\alpha)}$ onto $M$ such that
\begin{eqnarray}
\psi_{\alpha}(Y_{\alpha}(a,z)u)=Y(\Delta(\alpha,z)a,z)\psi_{\alpha}(u)\;\;\;
\mbox{for }a\in V, u\in M^{(\alpha)}.
\end{eqnarray}
For $\alpha=0$, we may choose $\psi_{0}={\rm id}_{M}$.

Set $U=\oplus _{\alpha\in L}V^{(\alpha)}$. For $\alpha,\beta\in L$, we
have a linear isomorphism $\psi_{\beta-\alpha}^{-1}\psi_{\beta}$ from
$V^{(\beta)}$ onto $V^{(\beta-\alpha)}$ satisfying the following condition:
\begin{eqnarray}\label{3.10}
\psi_{\beta-\alpha}^{-1}\psi_{\beta}(Y_{\b}(a,z)u)=Y_{\b-\a}(\Delta(\alpha,z)a,z)
\psi_{\beta-\alpha}^{-1}\psi_{\beta}(u)
\;\;\mbox{ for }a\in V, u\in V^{(\beta)}.
\end{eqnarray}
Then we may extend $\psi_{\alpha}$ to an automorphism of $U$ such that
\begin{eqnarray}\label{3.11}
\psi_{\alpha}(Y_{\beta}(a,z)u)=Y_{\beta+\a}(\Delta(\alpha,z)a,z)\psi_{\alpha}(u)
\;\;\;\mbox{ for }a\in V, u\in V^{(\a)}.
\end{eqnarray}
Then it is easy to see that $\psi_{\alpha+\beta}=\psi_{\alpha}\psi_{\beta}$
for any $\alpha, \beta\in L$. In other words, $\psi$ gives rise to a
representation of $L$ on $U.$

By a simple calculation we get:
\begin{equation}\label{3.12}
\Delta(\alpha,z)\beta=\beta+z^{-1}\langle \alpha,\beta\rangle{\bf 1},\;\;
\Delta(\alpha,z)\omega=\omega+z^{-1}\alpha+
z^{-2}\frac{\langle \alpha,\alpha\rangle}{2}{\bf 1}.
\end{equation}
\begin{lem}\label{l3.4}
{\it For any $\alpha\in L, h\in H$, we have}
\begin{eqnarray}
& &\psi_{\alpha}h(n)=h(n)\psi_{\alpha}+\delta_{n,0}\langle\alpha,h\rangle
\;\;\;\mbox{for }n\in {\Z},\label{3.13}\\
& &\psi_{\alpha}\Delta(h,z)=z^{\langle\alpha,h\rangle}\Delta(h,z)\psi_{\alpha}
,\label{3.14}\\
&
&\psi_{\alpha}e^{zL(-1)}\psi_{-\alpha}e^{-zL(-1)}=E^{-}(\alpha,z).\label{3.15}
\end{eqnarray}
\end{lem}

{\bf Proof.} By definition, we get
$\Delta(h,z)\alpha=\alpha +z^{-1}\langle\alpha,h\rangle {\bf 1}$. Then
(\ref{3.13}) is clear
and (\ref{3.14}) follows from (\ref{3.13}). By (\ref{3.11}) and (\ref{3.12}) we
get:
\begin{eqnarray}\label{3.16}
\psi_{\alpha}L(-1)=(L(-1)+\alpha(-1))\psi_{\alpha}.
\end{eqnarray}
Thus $\psi_{\alpha}e^{zL(-1)}\psi_{\alpha}^{-1}=e^{z(L(-1)+\alpha(-1))}$.
Then (\ref{3.15}) easily follows from Lemma 3.1.$\;\;\;\;\Box$

For any $\alpha\in L, h\in H$, we define:
\begin{eqnarray}
V^{(\alpha,h)}=\{u\in V^{(\alpha)}|h'(0)u=\langle h',h+\alpha\rangle u\;\;
\mbox{ for }h'\in H\}.
\end{eqnarray}
Then by Lemma \ref{l3.4} we have:
\begin{eqnarray}
V^{(\alpha)}=\oplus_{h\in H}V^{(\alpha,h)},\;
\psi_{\alpha}V^{(\alpha,h)}=V^{(0,h)}\;\;\mbox{ for }h\in H.
\end{eqnarray}

Let $P=\{\lambda\in H|V^{(0,\lambda)}\ne 0\}$. As $V$ is simple it
is easy to prove that $P$ is a subgroup of $H$ (cf. [LX]). Let $A=L\times P$
be the product group.
We define:
\be{3.19}
\eta((\alpha_{1},\lambda_{1}),(\alpha_{2},\lambda_{2}))=
-\langle \alpha_{1},\alpha_{2}\rangle
-\langle \alpha_{1},\lambda_{2}\rangle-\langle \alpha_{2},\lambda_{1}\rangle
\in {1\over T}{\Z}/2{\Z},
\ee
\be{3.20}
C((\alpha_{1},\lambda_{1}),(\alpha_{2},\lambda_{2}))=
e^{(\langle \alpha_{1},\lambda_{2}\rangle-\langle \alpha_{2},\lambda_{1}
\rangle)\pi i}\in {\C}^{*}
\ee
for any $(\alpha_{i},\lambda_{i})\in A, i=1,2$. Then $\eta(\cdot,\cdot)$ and
$C(\cdot,\cdot)$ satisfy the following conditions:
\be{3.21}
\eta(a,b)=\eta(b,a),\;\;\eta(a+b,c)=\eta(a,c)+\eta(b,c)
\ee
\be{3.22}
C(a,a)=1,\;\;C(a,b)=C(b,a)^{-1},\;\;C(a+b,c)=C(a,c)C(b,c)
\ee
for $a,b,c\in G$. $C(\cdot,\cdot)$ will be our commutator map later.

\bd{d3.5} For $u\in V^{(\alpha)}, v\in V^{(\beta)},
\alpha,\beta\in L$, we define $Y_{\a}(u,z)v\in V^{(\alpha+\beta)}\{z\}$ as
follows:
\be{3.23}
Y_{\a}(u,z)v=
\psi_{-\alpha-\beta}E^{-}(\alpha,z)Y(\psi_{\alpha}\Delta(\beta,z)u,z)
\Delta(\alpha,-z)\psi_{\beta}(v).
\ee
\ed

Set $U=\oplus_{{\a}\in L}V^{(\a)}.$ Then this defines a map $Y(\cdot,z)$ from
$U$ to $(\End(U))\{z\}$ via $Y(u,z)=Y_{\a}(u,z)$ for
$u\in V^{(\a)}.$ Notice that for any
$u\in V^{(\alpha, h_{1})}, v\in V^{(\beta,h_{2})}$,
$Y(u,z)v\in V^{(\alpha+\beta, h_{1}+h_{2})}\{z\}$.

\bp{p3.6} {\it The following $L(-1)$-derivative property holds:}
\begin{eqnarray}\label{3.24}
Y(L(-1)u,z)v={d\over dz}Y(u,z)v
\end{eqnarray}
for any $u,v\in U.$
\ep

{\bf Proof.} Let $\alpha,\beta\in L$ and let $u\in V^{(\alpha)},
v\in V^{(\beta)}$. Then
\begin{eqnarray}
& &Y(L(-1)u,z)v=
\psi_{-\alpha-\beta}E^{-}(\alpha,z)Y(\psi_{\alpha}\Delta(\beta,z)L(-1)u,z)
\Delta(\alpha,-z)\psi_{\beta}(v)\nonumber\\
& &\ \ \ \ =\psi_{-\alpha-\beta}E^{-}(\alpha,z)
Y(\psi_{\alpha}[\Delta(\beta,z),L(-1)]u,z)
\Delta(\alpha,-z)\psi_{\beta}(v)\nonumber\\
& &\ \ \ \ \ \ \ +\psi_{-\alpha-\beta}E^{-}(\alpha,z)
Y([\psi_{\alpha},L(-1)]\Delta(\beta,z)u,z)
\Delta(\alpha,-z)\psi_{\beta}(v)\nonumber\\
& &\ \ \ \ \ \ \
+\psi_{-\alpha-\beta}E^{-}(\alpha,z)Y(L(-1)\psi_{\alpha}\Delta(\beta,z)u,z)
\Delta(\alpha,-z)\psi_{\beta}(v).\label{3.25}
\end{eqnarray}
Note that (\ref{3.16}) is equivalent to
$$[\psi_{\alpha},L(-1)]=\alpha(-1)\psi_{\alpha}.$$
{}From (\ref{11.rjac}) with $u=\a(-1)\cdot {\bf 1},$ which is
$\sigma_h$-invariant
for any $h\in H,$ we have
\begin{eqnarray*}
& &\ \ \ Y(\alpha(-1)w,z)\nonumber\\
& &=\sum_{i=0}^{\infty}
\left(\begin{array}{c}-1\\i\end{array}\right)\left( (-z)^{i}\alpha(-1-i)Y(w,z)
+z^{-1-i}Y(w,z)\alpha(i)\right)\nonumber\\
& &=\sum_{i=0}^{\infty}\left(z^{i}\alpha(-1-i)Y(w,z)+
z^{-1-i}Y(w,z)\alpha(i)\right)\nonumber\\
& &=\alpha(z)^{-}Y(w,z)+Y(w,z)\alpha(z)^{+}\label{3.27}
\end{eqnarray*}
where
$$\alpha(z)^{-}=\sum_{n=0}^{\infty}\alpha(-n-1)z^n,
\alpha(z)^{+}=\sum_{n=1}^{\infty}\alpha(n)z^{-n-1}.$$
Thus
\begin{eqnarray*}
& &\ \ \psi_{-\alpha-\beta}E^{-}(\alpha,z)
Y([\psi_{\alpha},L(-1)]\Delta(\beta,z)u,z)
\Delta(\alpha,-z)\psi_{\beta}(v)\nonumber\\
& &=\psi_{-\alpha-\beta}E^{-}(\alpha,z)\alpha(z)^{-}Y(\psi_{\alpha}
\Delta(\beta,z)u,z)\Delta(\alpha,-z)\psi_{\beta}(v)\nonumber\\
& &\ \ +\psi_{-\alpha-\beta}E^{-}(\alpha,z)Y(\psi_{\alpha}
\Delta(\beta,z)u,z)\alpha(z)^{+}\Delta(\alpha,-z)\psi_{\beta}(v)\nonumber\\
& &=\psi_{-\alpha-\beta}\left(\frac{d}{dz}E^{-}(\alpha,z)\right)Y(\psi_{\alpha}
\Delta(\beta,z)u,z)\Delta(\alpha,-z)\psi_{\beta}(v)\nonumber\\
& &\ \ +\psi_{-\alpha-\beta}E^{-}(\alpha,z)Y(\psi_{\alpha}
\left(\frac{d}{dz}\Delta(\beta,z)u,z)\right)\Delta(\alpha,-z)\psi_{\beta}(v).
\label{3.28}
\end{eqnarray*}
Now from (\ref{li3}) and Proposition \ref{p2.15}
we obtain
\begin{eqnarray*}
Y(L(-1)u,z)v={d\over dz}Y(u,z).\;\;\;\;\Box
\end{eqnarray*}

\bt{t3.7} {\it For any} $u\in V^{(\alpha,h_{1})},
v\in V^{(\beta,h_{2})}, w\in V^{(\gamma,h_{3})},
\alpha,\beta,\gamma\in L, h_{1},h_{2},h_{3}\in P$, {\it we have the following
generalized Jacobi identity:}
\begin{eqnarray}
& &z_{0}^{-1}\delta\left(\frac{z_{1}-z_{2}}{z_{0}}\right)
\left(\frac{z_{1}-z_{2}}{z_{0}}\right)^{\eta((\alpha,h_{1}),(\beta,h_{2}))}
Y(u,z_{1})Y(v,z_{2})w\nonumber\\
&-&C((\alpha,h_{1}),(\beta,h_{2}))z_{0}^{-1}\delta\left(\frac{z_{2}-z_{1}}
{-z_{0}}\right)
\left(\frac{z_{2}-z_{1}}{z_{0}}\right)^{\eta((\alpha,h_{1}),(\beta,h_{2}))}
Y(v,z_{2})Y(u,z_{1})w\nonumber\\
&=&z_{2}^{-1}\delta\left(\frac{z_{1}-z_{0}}{z_{2}}\right)
\left(\frac{z_{2}+z_{0}}{z_{1}}\right)^{\eta((\alpha,h_{1}),(\gamma,h_{3}))}
Y(Y(u,z_{0})v,z_{2})w.\label{3.31}
\end{eqnarray}
Moreover, $(U, {\bf 1},\omega, Y, T, A, \eta(\cdot,\cdot),
C(\cdot,\cdot))$ is a
generalized vertex algebra in the sense of [DL].
\et

{\bf Proof.} By (\ref{3.21}), (\ref{3.22}) and Proposition \ref{p3.6} all
the axioms in the definition of a generalized vertex algebra (Chapter 9
of [DL]) hold except the Jacobi identity.
By Definition \ref{d3.5} and Lemma \ref{l3.1} together with the relation
(\ref{3.14})
we have:
\begin{eqnarray}
& &\ \ \ Y(u,z_{1})Y(v,z_{2})w\nonumber\\
& &=z_{1}^{\langle \alpha,\beta+\gamma\rangle}
\psi_{-\alpha-\beta-\gamma}E^{-}(\alpha,z_{1})
Y\left(\Delta(\beta+\gamma,z_{1})\psi_{\alpha}u,z_{1}\right)
\Delta(\alpha,-z_{1})
\psi_{\beta+\gamma}Y(v,z_{2})w\nonumber\\
& &=z_{1}^{\langle \alpha,\beta+\gamma\rangle}
z_{2}^{\langle \beta,\gamma\rangle}
\psi_{-\alpha-\beta}E^{-}(\alpha,z_{1})
Y\left(\Delta(\beta+\gamma,z_{1})\psi_{\alpha}u,z_{1}\right)
\Delta(\alpha,-z_{1})\nonumber\\
& &\ \ \ \cdot E^{-}(\beta,z_{2})Y(\Delta(\gamma,z_{2})
\psi_{\beta}v,z_{2})\Delta(\beta,-z_{2})\psi_{\gamma}w\nonumber\\
& &=\left(1-{z_{2}\over z_{1}}\right)^{\langle \alpha,\beta\rangle}
z_{1}^{\langle \alpha,\beta+\gamma\rangle}z_{2}^{\langle \beta,\gamma\rangle}
\psi_{-\alpha-\beta-\gamma}\nonumber\\
& &\ \ \ \cdot E^{-}(\alpha,z_{1})E^{-}(\beta,z_{2})Y(\Delta(\beta,z_{1}-z_{2})
\Delta(\gamma,z_{1})\psi_{\alpha}u,z_{1})\nonumber\\
& &\ \  \ \cdot Y\left(\Delta(\alpha,-z_{1}+z_{2})
\Delta(\gamma,z_{2}\right)\psi_{\beta}v,z_{2})\Delta(\alpha,-z_{1})
\Delta(\beta,-z_{2})\psi_{\gamma}w\nonumber\\
& &=
(z_{1}-z_{2})^{\langle \alpha,\beta\rangle}z_{1}^{\langle \alpha,\gamma\rangle}
z_{2}^{\langle \beta,\gamma\rangle}
\psi_{-\alpha-\beta-\gamma}\nonumber\\
& &\ \ \ \cdot
E^{-}(\alpha,z_{1})E^{-}(\beta,z_{2})Y\left(\Delta(\beta,z_{1}-z_{2})
\Delta(\gamma,z_{1})\psi_{\alpha}u,z_{1}\right)\nonumber\\
& &\ \ \ \cdot Y\left(\Delta(\alpha,-z_{1}+z_{2})
\Delta(\gamma,z_{2})\psi_{\beta}v,z_{2}\right)\Delta(\alpha,-z_{1})
\Delta(\beta,-z_{2})\psi_{\gamma}w.
\end{eqnarray}
Symmetrically,
\begin{eqnarray}
& & \ \ \ Y(v,z_{2})Y(u,z_{1})w\nonumber\\
& &=
(z_{2}-z_{1})^{\langle \alpha,\beta\rangle}z_{1}^{\langle \alpha,\gamma\rangle}
z_{2}^{\langle \beta,\gamma\rangle}
\psi_{-\alpha-\beta-\gamma}\nonumber\\
& &\ \ \ \cdot
E^{-}(\alpha,z_{1})E^{-}(\beta,z_{2})Y\left(\Delta(\alpha,z_{2}-z_{1})
\Delta(\gamma,z_{2})\psi_{\beta}v,z_{2}\right)\nonumber\\
& &\ \ \ \cdot Y\left(\Delta(\beta,-z_{2}+z_{1})
\Delta(\gamma,z_{1})\psi_{\alpha}u,z_{1}\right)\Delta(\alpha,-z_{1})
\Delta(\beta,-z_{2})\psi_{\gamma}w.
\end{eqnarray}
Since
$$\beta(0)\Delta(\gamma,z_{1})\psi_{\alpha}u=\Delta(\gamma,z_{1})\beta(0)\psi_{\alpha}u
=\langle \beta,h_{1}\rangle\Delta(\gamma,z_{1})\psi_{\alpha}u,$$
 we see that
$(z_{1}-z_{2})^{-\langle \beta,h_{1}\rangle}\Delta(\beta,z_{1}-z_{2})
\Delta(\gamma,z_{1})\psi_{\alpha}u$ involves only integral powers of
$(z_{1}-z_{2})$.
Similarly, $(z_{1}-z_{2})^{-\langle
\alpha,h_{2}\rangle}\Delta(\alpha,-z_{1}+z_{2})
\Delta(\gamma,z_{2})\psi_{\beta}v$ involves only integral powers of
$(z_{1}-z_{2})$.
Using properties of $\delta$-functions (cf. [FLM]) we see that
\begin{eqnarray}
& &\ \ \ z_{0}^{-1}\delta\left(\frac{z_{1}-z_{2}}{z_{0}}\right)
Y(u,z_{1})Y(v,z_{2})w\nonumber\\
& &=z_{0}^{-1}\delta\left(\frac{z_{1}-z_{2}}{z_{0}}\right)
(z_{1}-z_{2})^{\langle \alpha,\beta\rangle}z_{1}^{\langle \alpha,\gamma\rangle}
z_{2}^{\langle \beta,\gamma\rangle}
\psi_{-\alpha-\beta-\gamma}\nonumber\\
& &\ \ \ \cdot E^{-}(\alpha,z_{1})E^{-}(\beta,z_{2})
\left(\frac{z_{1}-z_{2}}{z_{0}}\right)^{\langle \beta,h_{1}\rangle}
Y\left(\Delta(\beta,z_{0})
\Delta(\gamma,z_{1})\psi_{\alpha}u,z_{1}\right)\nonumber\\
& &\ \ \ \cdot \left(\frac{z_{1}-z_{2}}{z_{0}}\right)^{\langle
\alpha,h_{2}\rangle}
Y\left(\Delta(\alpha,-z_{0})
\Delta(\gamma,z_{2})\psi_{\beta}v,z_{2}\right)\Delta(\alpha,-z_{1})
\Delta(\beta,-z_{2})\psi_{\gamma}w\nonumber\\
& &=z_{0}^{-1}\delta\left(\frac{z_{1}-z_{2}}{z_{0}}\right)
\left(\frac{z_{1}-z_{2}}{z_{0}}\right)^{\langle \alpha,\beta\rangle+
\langle \alpha,h_{2}\rangle+\langle \beta,h_{1}\rangle}
z_{1}^{\langle \alpha,\gamma\rangle}z_{2}^{\langle \beta,\gamma\rangle}
z_{0}^{\langle \alpha,\beta\rangle}\psi_{-\alpha-\beta-\gamma}\nonumber\\
& &\ \ \ \cdot E^{-}(\alpha,z_{1})E^{-}(\beta,z_{2})Y\left(\Delta(\beta,z_{0})
\Delta(\gamma,z_{1})\psi_{\alpha}u,z_{1}\right)\nonumber\\
& &\ \ \ \cdot Y\left(\Delta(\alpha,-z_{0})
\Delta(\gamma,z_{2})\psi_{\beta}v,z_{2}\right)\Delta(\alpha,-z_{1})
\Delta(\beta,-z_{2})\psi_{\gamma}w,
\end{eqnarray}
and that
\begin{eqnarray}
& &\ \ \ z_{0}^{-1}\delta\left(\frac{z_{2}-z_{1}}{-z_{0}}\right)
Y(v,z_{2})Y(u,z_{1})w\nonumber\\
& &=z_{0}^{-1}\delta\left(\frac{z_{2}-z_{1}}{-z_{0}}\right)
(z_{2}-z_{1})^{\langle \alpha,\beta\rangle}z_{1}^{\langle \alpha,\gamma\rangle}
z_{2}^{\langle \beta,\gamma\rangle}
\psi_{-\alpha-\beta-\gamma}\nonumber\\
& &\ \ \ \cdot
E^{-}(\alpha,z_{1})E^{-}(\beta,z_{2})Y\left(\Delta(\alpha,z_{2}-z_{1})
\Delta(\gamma,z_{2})\psi_{\beta}v,z_{2}\right)\nonumber\\
& &\ \ \ \cdot Y\left(\Delta(\beta,-z_{2}+z_{1})
\Delta(\gamma,z_{1})\psi_{\alpha}u,z_{1}\right)\Delta(\alpha,-z_{1})
\Delta(\beta,-z_{2})\psi_{\gamma}w\nonumber\\
& &=z_{0}^{-1}\delta\left(\frac{z_{2}-z_{1}}{-z_{0}}\right)
(z_{2}-z_{1})^{\langle \alpha,\beta\rangle}z_{1}^{\langle \alpha,\gamma\rangle}
z_{2}^{\langle \beta,\gamma\rangle}
\psi_{-\alpha-\beta-\gamma}E^{-}(\alpha,z_{1})E^{-}(\beta,z_{2})\nonumber\\
& &\ \ \ \cdot \left(\frac{z_{2}-z_{1}}{-z_{0}}\right)^{\langle
\alpha,h_{2}\rangle}
Y\left(\Delta(\alpha,-z_{0})
\Delta(\gamma,z_{2})\psi_{\beta}v,z_{2}\right)\nonumber\\
& &\ \ \ \cdot \left(\frac{-z_{2}+z_{1}}{z_{0}}\right)^{\langle
\beta,h_{1}\rangle}
Y\left(\Delta(\beta,z_{0})
\Delta(\gamma,z_{1})\psi_{\alpha}u,z_{1}\right)\Delta(\alpha,-z_{1})
\Delta(\beta,-z_{2})\psi_{\gamma}w\nonumber\\
& &=e^{(\langle \beta,h_{1}\rangle-\langle \alpha,h_{2}\rangle)\pi i}
z_{0}^{-1}\delta\left(\frac{z_{2}-z_{1}}{-z_{0}}\right)
\left(\frac{z_{2}-z_{1}}{z_{0}}\right)^{\langle \alpha,\beta\rangle+
\langle \alpha,h_{2}\rangle+\langle \beta,h_{1}\rangle}
z_{1}^{\langle \alpha,\gamma\rangle}z_{2}^{\langle \beta,\gamma\rangle}
z_{0}^{\langle \alpha,\beta\rangle}\psi_{-\alpha-\beta-\gamma}
\nonumber\\
& &\ \ \ \cdot
E^{-}(\alpha,z_{1})E^{-}(\beta,z_{2})Y\left(\Delta(\alpha,-z_{0})
\Delta(\gamma,z_{2})\psi_{\beta}v,z_{2}\right)\nonumber\\
& &\ \ \ \cdot Y\left(\Delta(\beta,z_{0})
\Delta(\gamma,z_{1})\psi_{\alpha}u,z_{1}\right)\Delta(\alpha,-z_{1})
\Delta(\beta,-z_{2})\psi_{\gamma}w.
\end{eqnarray}
Thus
\begin{eqnarray}
& &\ \ \ z_{0}^{-1}\delta\left(\frac{z_{1}-z_{2}}{z_{0}}\right)
\left(\frac{z_{1}-z_{2}}{z_{0}}\right)^{\eta((\alpha,h_{1}),(\beta,h_{2}))}
Y(u,z_{1})Y(v,z_{2})w\nonumber\\
& &\ \ \
-C((\alpha,h_{1}),(\beta,h_{2}))z_{0}^{-1}\delta\left(\frac{z_{2}-z_{1}}{-z_{0}}\right)
\left(\frac{z_{2}-z_{1}}{z_{0}}\right)^{\eta((\alpha,h_{1}),(\beta,h_{2}))}
Y(v,z_{2})Y(u,z_{1})w\nonumber\\
& &=z_{2}^{-1}\delta\left(\frac{z_{1}-z_{0}}{z_{2}}\right)
z_{1}^{\langle \alpha,\gamma\rangle}z_{2}^{\langle \beta,\gamma\rangle}
z_{0}^{-\langle \alpha, \beta\rangle}\psi_{-\alpha-\beta-\gamma}
E^{-}(\alpha,z_{1})E^{-}(\beta,z_{2})\nonumber\\
& &\ \ \ \cdot
Y(Y(\Delta(\beta,z_{0})\Delta(\gamma,z_{1})\psi_{\alpha}(u),z_{0})
\Delta(\alpha,-z_{0})\Delta(\gamma,z_{2})\psi_{\beta}(v),z_{2})\nonumber\\
& &\ \ \ \cdot \Delta(\alpha,-z_{1})\Delta(\beta,-z_{2})\psi_{\gamma}(w).
\end{eqnarray}
On the other hand,
\begin{eqnarray}
& &\ \ \ Y(Y(u,z_{0})v,z_{2})w\nonumber\\
& &=E^{-}(\alpha+\beta,z_{2})\psi_{-\alpha-\beta-\gamma}
Y\left(\psi_{\alpha+\beta}\Delta(\gamma,z_{2})Y(u,z_{0})v,z_{2}\right)
\Delta(\alpha+\beta,-z_{2})\psi_{\gamma}w\nonumber\\
& &=E^{-}(\alpha+\beta,z_{2})\psi_{-\alpha-\beta-\gamma}
\nonumber\\
& &\ \ \ \cdot
Y\left(\psi_{\alpha+\beta}\Delta(\gamma,z_{2})E^{-}(\alpha,z_{0})
\psi_{-\alpha-\beta}Y(\psi_{\alpha}\Delta(\beta,z_{0})
u,z_{0})\Delta(\alpha,,-z_{0})\psi_{\beta}v,z_{2}\right)\nonumber\\
& &\ \ \ \cdot \Delta(\alpha+\beta,-z_{2})\psi_{\gamma}w\nonumber\\
& &=z_{2}^{\langle \alpha+\beta,\gamma\rangle}
z_{0}^{\langle \alpha,\beta\rangle}
E^{-}(\alpha+\beta,z_{2})\psi_{-\alpha-\beta-\gamma}\nonumber\\
& &\ \ \ \cdot
Y\left(\Delta(\gamma,z_{2})E^{-}(\alpha,z_{0})Y(\Delta(\beta,z_{0})
\psi_{\alpha}u,z_{0})\Delta(\alpha,-z_{0})\psi_{\beta}v,z_{2}\right)
\nonumber\\
& &\ \ \ \cdot\Delta(\alpha+\beta,-z_{2})\psi_{\gamma}w\nonumber\\
& &=\left(1+{z_{0}\over z_{2}}\right)^{\langle \alpha,\gamma\rangle}
z_{2}^{\langle \alpha+\beta,\gamma\rangle}z_{0}^{\langle \alpha,\beta\rangle}
E^{-}(\alpha+\beta,z_{2})\psi_{-\alpha-\beta-\gamma}\nonumber\\
& &\ \ \ \cdot
Y\left(E^{-}(\alpha,z_{0})\Delta(\gamma,z_{2})Y(\Delta(\beta,z_{0})
\psi_{\alpha}u,z_{0})\Delta(\alpha,-z_{0})\psi_{\beta}v,z_{2}\right)
\nonumber\\
& &\ \ \ \cdot\Delta(\alpha+\beta,-z_{2})\psi_{\gamma}w\nonumber\\
& &=(z_{2}+z_{0})^{\langle \alpha,\gamma\rangle}
z_{2}^{\langle \beta,\gamma\rangle}z_{0}^{\langle \alpha,\beta\rangle}
E^{-}(\alpha+\beta,z_{2})E^{-}(\alpha,z_{0}+z_{2})E^{-}(-\alpha,z_{2})
\psi_{-\alpha-\beta-\gamma}\nonumber\\
& &\ \ \ \cdot Y(\Delta(\gamma,z_{2})Y(\Delta(\beta,z_{0})
\psi_{\alpha}u,z_{0})\Delta(\alpha,-z_{0})\psi_{\beta}v,z_{2})
\nonumber\\
& &\ \ \ \cdot z_{2}^{-\alpha(0)}(z_{2}+z_{0})^{\alpha (0)}E^{+}(\alpha,z_{2})
E^{+}(-\alpha,z_{2}+z_{0})
\Delta(\alpha+\beta,-z_{2})\psi_{\gamma}w\nonumber\\
& &=(z_{2}+z_{0})^{\langle \alpha,\gamma\rangle}
z_{2}^{\langle \beta,\gamma\rangle}z_{0}^{\langle \alpha,\beta\rangle}
E^{-}(\beta,z_{2})E^{-}(\alpha,z_{0}+z_{2})\psi_{-\alpha-\beta-\gamma}\nonumber\\
& &\ \ \ \cdot Y\left(\Delta(\gamma,z_{2})Y(\Delta(\beta,z_{0})
\psi_{\alpha}u,z_{0})\Delta(\alpha,-z_{0})\psi_{\beta}v,z_{2}\right)
\nonumber\\
& &\ \ \ \cdot z_{2}^{-\langle \alpha, h_{3}\rangle}
(z_{2}+z_{0})^{\langle \alpha, h_{3}\rangle}E^{+}(-\beta,z_{2})
E^{+}(-\alpha,z_{2}+z_{0})(-z_{2})^{\langle \alpha+\beta, h_{3}\rangle}
\psi_{\gamma}w\nonumber\\
& &=(z_{2}+z_{0})^{\langle \alpha,\gamma\rangle}
z_{2}^{\langle \beta,\gamma\rangle}z_{0}^{\langle \alpha,\beta\rangle}
E^{-}(\beta,z_{2})E^{-}(\alpha,z_{0}+z_{2})\psi_{-\alpha-\beta-\gamma}\nonumber\\
& &\ \ \ \cdot Y\left(Y(\Delta(\gamma,z_{2}+z_{0})\Delta(\beta,z_{0})
\psi_{\alpha}u,z_{0})\Delta(\gamma,z_{2})\Delta(\alpha,-z_{0})\psi_{\beta}v,z_{2}\right)
\nonumber\\
& &\ \ \ \cdot z_{2}^{-\langle \alpha, h_{3}\rangle}
(z_{2}+z_{0})^{\langle \alpha, h_{3}\rangle}E^{+}(-\beta,z_{2})
E^{+}(-\alpha,z_{2}+z_{0})(-z_{2})^{\langle \alpha+\beta, h_{3}\rangle}
\psi_{\gamma}w.
\end{eqnarray}
Hence
\begin{eqnarray}
& &\ \ \
z_{2}^{-1}\delta\left(\frac{z_{1}-z_{0}}{z_{2}}\right)Y(Y(u,z_{0})v,z_{2})w
\nonumber\\
& &=z_{2}^{-1}\delta\left(\frac{z_{1}-z_{0}}{z_{2}}\right)
(z_{2}+z_{0})^{\langle \alpha,\gamma\rangle}
z_{2}^{\langle \beta,\gamma\rangle}z_{0}^{\langle \alpha,\beta\rangle}
E^{-}(\beta,z_{2})E^{-}(\alpha,z_{1})\psi_{-\alpha-\beta-\gamma}\nonumber\\
& &\ \ \ \cdot Y\left(Y(\Delta(\gamma,z_{2}+z_{0})\Delta(\beta,z_{0})
\psi_{\alpha}u,z_{0})\Delta(\gamma,z_{2})\Delta(\alpha,-z_{0})
\psi_{\beta}v,z_{2}\right)
\nonumber\\
& &\ \ \ \cdot z_{2}^{-\langle \alpha, h_{3}\rangle}
(z_{2}+z_{0})^{\langle \alpha, h_{3}\rangle}E^{+}(-\beta,z_{2})
E^{+}(-\alpha,z_{1})(-z_{2})^{\langle \alpha+\beta, h_{3}\rangle}
\psi_{\gamma}w\nonumber\\
& &=z_{2}^{-1}\delta\left(\frac{z_{1}-z_{0}}{z_{2}}\right)
(z_{2}+z_{0})^{\langle \alpha,\gamma\rangle}
z_{2}^{\langle \beta,\gamma\rangle}z_{0}^{\langle \alpha,\beta\rangle}
E^{-}(\beta,z_{2})E^{-}(\alpha,z_{1})\psi_{-\alpha-\beta-\gamma}\nonumber\\
& &\ \ \ \cdot \left(\frac{z_{2}+z_{0}}{z_{1}}\right)^{\langle
\gamma,h_{1}\rangle}
Y\left(Y(\Delta(\gamma,z_{1})\Delta(\beta,z_{0})
\psi_{\alpha}u,z_{0})\Delta(\gamma,z_{2})
\Delta(\alpha,-z_{0})\psi_{\beta}v,z_{2}\right)
\nonumber\\
& &\ \ \ \cdot \left(\frac{z_{2}+z_{0}}{z_{1}}\right)^{\langle
\alpha,h_{3}\rangle}
\Delta(\alpha,-z_{1})
\Delta(\beta,-z_{2})\psi_{\gamma}w\nonumber\\
& &=z_{2}^{-1}\delta\left(\frac{z_{1}-z_{0}}{z_{2}}\right)
\left(\frac{z_{2}+z_{0}}{z_{1}}\right)^{-\eta((\alpha,h_{1}),(\gamma,h_{3}))}
z_{1}^{\langle \alpha,\gamma\rangle}z_{2}^{\langle \beta,\gamma\rangle}
z_{0}^{\langle \alpha,\beta\rangle}E^{-}(\beta,z_{2})E^{-}(\alpha,z_{1})
\nonumber\\
& &\ \ \ \cdot \psi_{-\alpha-\beta-\gamma}Y\left(Y(\Delta(\gamma,z_{1})
\Delta(\beta,z_{0})
\psi_{\alpha}u,z_{0})\Delta(\gamma,z_{2})
\Delta(\alpha,-z_{0})\psi_{\beta}v,z_{2}\right)
\nonumber\\
& &\ \ \ \cdot \Delta(\alpha,-z_{1})\Delta(\beta,-z_{2})\psi_{\gamma}w.
\end{eqnarray}
The generalized Jacobi identity (\ref{3.31}) now follows immediately.
$\;\;\;\;\Box$

Next we shall extend certain $V$-modules to modules for $U$ considered
as a generalized vertex algebra via Theorem 3.5.
Let $M$ be an irreducible $V$-module (with finite-dimensional homogeneous
subspaces).
Since $[L(0), h(0)]=0$ for any $h\in H$, $H$ preserves each homogeneous
subspace of $M$
so that there exist $0\ne u\in M, \lambda\in H^{*}$ such that
$h(0)u=\lambda(h)u$ for $h\in H$. Since $H$ acts semisimply on $V$ (by
assumption) and
$u$ generates $M$ by $V$ (from the irreducibility of $M$), $H$ also acts
semisimply on $M$.
For any $\lambda\in H^{*}$, we define
\begin{eqnarray}
M^{(0,\lambda)}=\{u\in M|h(0)u=\lambda(h)u\;\;\mbox{for }h\in H\}.
\end{eqnarray}

Set
\begin{eqnarray}
P(M)=\{\lambda\in H^{*}|M^{(0,\lambda)}\ne 0\}.
\end{eqnarray}
Since $M$ is irreducible, $P(M)$ is an irreducible $P(V)$-set. Thus $L\times
P(M)$ is an
irreducible $(L\times P(V))$-set.
{\em Suppose that for any $\alpha\in L$, $\alpha(0)$ has rational eigenvalues
on $M$.}
Let $\lambda_{0}\in P(M)$ be any $H$-weight of $M$. Then
$P(M)=\lambda_{0}+P(V)$. By using a basis of $L$,
we see that there is a positive integer $K$ such that
$\langle \lambda_{0},\alpha\rangle\in {1\over K}{\Z}$ for any $\alpha\in L$.
Therefore
\begin{eqnarray}
\langle \lambda,\alpha\rangle\in {1\over TK}{\Z}\;\;\;\mbox{for any }
\lambda\in P(M),\alpha\in L.
\end{eqnarray}
Using formula (3.19) we extend the definition of
$\eta(\cdot,\cdot)$ to $(L\times P(M))\times (L\times P(M))$ with values
in ${1\over TK}\Bbb{Z}$.

Recall that $(M^{(\alpha)},Y_{\alpha}(\cdot,z))$ is a weak
$\sigma_{\alpha}$-twisted
$V$-module for any $\alpha\in L$.
Set $W=\oplus _{\alpha\in L}M^{(\alpha)}$.
For $a\in V^{(\alpha)}, u\in M^{(\beta)}, \alpha,\beta\in L$,
we define $Y_{W}(a,z)u\in M^{(\alpha+\beta)}\{z\}$ as follows:
\begin{eqnarray}
Y_{W}(a,z)u=\psi_{-\alpha-\beta}E^{-}(\alpha,z)Y_{M}(\psi_{\alpha}\Delta(\beta,z)a,z)
\Delta(\alpha,-z)\psi_{\beta}(u).
\end{eqnarray}
Then the same argument used in the proof of Theorem 3.5 shows that for any
$a\in V^{(\alpha,h_{1})},b\in V^{(\beta,h_{2})},u\in M^{(\gamma,h_{3})}$,
where $\alpha,\beta,\gamma\in L,h_{1},h_{2}\in P,h_{3}\in P(M)$, we have:
\begin{eqnarray}
& &z_{0}^{-1}\delta\left(\frac{z_{1}-z_{2}}{z_{0}}\right)
\left(\frac{z_{1}-z_{2}}{z_{0}}\right)^{\eta((\alpha,h_{1}),(\beta,h_{2}))}
Y_{W}(a,z_{1})Y_{W}(b,z_{2})u\nonumber\\
&-&C((\alpha,h_{1}),(\beta,h_{2}))z_{0}^{-1}\delta\left(\frac{z_{2}-z_{1}}
{-z_{0}}\right)
\left(\frac{z_{2}-z_{1}}{z_{0}}\right)^{\eta((\alpha,h_{1}),(\beta,h_{2}))}
Y_{W}(b,z_{2})Y_{W}(a,z_{1})u\nonumber\\
& &=z_{2}^{-1}\delta\left(\frac{z_{1}-z_{0}}{z_{2}}\right)
\left(\frac{z_{2}+z_{0}}{z_{1}}\right)^{\eta((\alpha,h_{1}),(\gamma,h_{3}))}
Y_{W}(Y(a,z_{0})b,z_{2})u.
\end{eqnarray}
Then we have:

\bt{t3.6}
{\it $(W,Y_{W})$ is a module for the generalized vertex  algebra
$U$ in the sense of [DL].}
\et

{}In view of Proposition 2.15, it is possible that various
 $V^{(\alpha)}$ in $U$ may be $V$-isomorphic to
each other. Next we shall reduce $U$ to a smaller space $\bar{U}$
such that the  multiplicity of any $\sigma_{\alpha}$-twisted
$V$-module $V^{(\alpha)}$ ($\alpha\in L$) is one,
and $\bar{U}$ is an abelian intertwining algebra in the sense of [DL]
rather than a generalized vertex algebra.

Set
\begin{eqnarray}
L_{0}=\{\alpha \in L| \sigma_{\alpha}={\rm id}_{V},\;V^{(\alpha)}\simeq V\}.
\end{eqnarray}

\begin{lem}\label{l3.7} {\it Let $\alpha,\beta\in L$. Then
$\sigma_{\alpha}=\sigma_{\beta}$
and $V^{(\alpha)}$ is isomorphic to
$V^{(\beta)}$ as a weak $\sigma_{\alpha}$-twisted $V$-module if and only if
$\alpha-\beta\in L_{0}$.}
\end{lem}

{\bf Proof.} Suppose that $\sigma_{\alpha}=\sigma_{\beta}$
and $V^{(\alpha)}$ is isomorphic to
$V^{(\beta)}$ as a weak $\sigma_{\alpha}$-twisted $V$-module. Let
$\phi$ be a $V$-isomorphism from $V^{(\alpha)}$
onto $V^{(\beta)}$.
Then $\psi_{\beta}\phi \psi_{-\beta}$ is a linear isomorphism from
$V^{(\alpha-\beta)}$ onto $V$ such that
\begin{eqnarray}
\psi_{\beta}\phi \psi_{-\beta}(Y(a,z)u)=Y(\Delta(\beta,z)\Delta(-\beta,z)a,z)
\psi_{\beta}\phi \psi_{-\beta}(u)=Y(a,z)\psi_{\beta}\phi \psi_{-\beta}(u)
\end{eqnarray}
for any $a\in V, u\in V^{(\alpha-\beta)}$. Then by definition,
$\alpha-\beta\in L_{0}$.

On the other hand, suppose that $\alpha-\beta\in L_{0}$ for some
$\alpha,\beta\in L$. Since
$\sigma_{\alpha-\beta}={\rm id}_{V}$, we have $\sigma_{\alpha}=\sigma_{\beta}$.
Let
$\psi$ be a $V$-isomorphism from $V^{(\alpha-\beta)}$ onto $V$. Then
$\psi_{-\beta}\psi \psi_{\beta}$ is a linear isomorphism
{}from $V^{(\alpha)}$ onto $V^{(\beta)}$ satisfying the condition:
\begin{eqnarray}
& &\;\;\;\;\psi_{\beta}^{-1}\psi \psi_{\beta}(Y_{\alpha}(a,z)u)
\nonumber\\
& &=\psi_{-\beta}\psi Y(\Delta(\beta,z)a,z)\psi_{\beta}u\nonumber\\
& &=\psi_{\beta}Y(\Delta(\beta,z)a,z)\psi\psi_{\beta}u\nonumber\\
&
&=Y_{\beta}(\Delta(-\beta,z)\Delta(\beta,z)a,z)\psi_{-\beta}\psi\psi_{\beta}u\nonumber\\
& &=Y_{\beta}(a,z)\psi_{-\beta}\psi\psi_{\beta}u
\end{eqnarray}
for any $a\in V$. Then $V^{(\alpha)}$ is $V$-isomorphic to $V^{(\beta)}$.
The proof is complete. $\;\;\;\;\Box$

For any $\alpha,\beta\in L_{0}$, by definition we have
$V^{(\alpha)}\simeq V\simeq V^{(\beta)}$. From Lemma 3.7,
$\alpha-\beta\in L_{0}$. Thus $L_{0}$ is a sublattice of $L$.

\begin{lem}\label{l3.8}
{\it $L_{0}$ is an even sublattice such that $L_{0}\subseteq P\cap P^{0}$,
where $P^{0}$
is the dual lattice of $P$.}
\end{lem}

{\bf Proof.} Let $\alpha\in L_{0}.$ Then
$e^{2\pi i \alpha(0)}={\rm id}_{V}$ if and only if $\alpha(0)$ has
integral eigenvalues on $V$. This proves $\alpha\in P^{0}$. From (\ref{3.13}),
we obtain:
\begin{eqnarray}
& &\psi_{\alpha}\beta(0)=\left(\beta(0)+\langle
\alpha,\beta\rangle\right)\psi_{\alpha},\\
& &\psi_{\alpha}L(0)=\left(L(0)+\alpha(0)+{1\over 2}\langle
\alpha,\alpha\rangle\right)
\psi_{\alpha}.
\end{eqnarray}
Then (3.42) implies that
$P(V^{(\alpha)})=-\alpha+P$. Since $V^{(\alpha)}$ is $V$-isomorphic to $V$,
$P=P(V^{(\alpha)})$. Thus $\alpha\in P$.
Let $u\in V^{(0,\lambda)}$ for $\lambda\in P$.
Then from (3.43) the $L(0)$-weight of $\psi_{-\alpha}(u)\in V^{(\alpha)}$ is
${\rm wt}u +\langle \alpha,\lambda\rangle+{1\over 2}\langle
\alpha,\alpha\rangle$.
In particular,
the $L(0)$-weight of $\psi_{-\alpha}({\bf 1})\in V^{(\alpha)}$ is
${1\over 2}\langle \alpha,\alpha\rangle$. Thus $\langle \alpha,\alpha\rangle\in
2{\Z}$. The proof is complete.$\;\;\;\;\Box$

Let $\alpha\in L_{0}$ and let
$\bar{\pi}_{\alpha}$ be a fixed $V$-isomorphism from $V^{(\alpha)}$ onto $V$.
If $\alpha=0$, we choose $\bar{\pi}_{0}={\rm id}_{V}$.
For any $\beta\in L$, considering the following composition map:
\begin{eqnarray}
V^{(\beta)}\longrightarrow V^{(\alpha)}\longrightarrow V\longrightarrow
V^{(\beta-\alpha)},
\end{eqnarray}
we obtain a linear isomorphism
$f=\psi_{\alpha-\beta}\bar{\pi}_{\alpha}\psi_{\beta-\alpha}$ from
$V^{(\beta)}$ onto $V^{(\beta-\alpha)}$. Then
\begin{eqnarray}
f(Y(a,z)u)=Y(\Delta(\alpha-\beta,z)\Delta(\beta-\alpha,z)a,z)f(u)=Y(a,z)f(u)
\end{eqnarray}
for any $a\in V, u\in V^{(\beta)}$.
Then we define an endomorphism $f_{\alpha}$ of $U$ as
follows:
\begin{eqnarray}
f_{\alpha}u= e^{\langle \alpha,\beta-\alpha\rangle\pi i}\psi_{\alpha-\beta}
\bar{\pi}_{\alpha}\psi_{\beta-\alpha}(u)\;\;\;\mbox{ for }
u\in V^{(\beta)}\subseteq U.
\end{eqnarray}
Since $\psi_{0}={\rm id}_{V}$, we have
 $f_{\alpha}|_{V^{(\alpha)}}=\bar{\pi}_{\alpha}$. Then we may use
$\bar{\pi}_{\alpha}$ to denote the linear isomorphism $f_{\alpha}$ of $U$
without any confusion. It is easy to see that $\bar{\pi}_{\alpha}$ satisfies
following condition:
\begin{eqnarray}
\bar{\pi}_{\alpha}(Y(a,z)u)=Y(a,z)\bar{\pi}_{\alpha}(u)\;\;\;
\mbox{ for }a\in V, u\in U.
\end{eqnarray}

\begin{lem}\label{l3.9} {\it For any $\alpha\in L_{0}, \beta\in L$, we have}
\begin{eqnarray}
\psi_{\beta}\bar{\pi}_{\alpha}
=e^{\langle \alpha,\beta\rangle \pi i}\bar{\pi}_{\alpha}\psi_{\beta}.
\end{eqnarray}
\end{lem}

{\bf Proof.} Let $\gamma\in L, u\in V^{(\gamma)}$. Then by (3.46) we have:
\begin{eqnarray}
\psi_{\beta}\bar{\pi}_{\alpha}(u)
&=&\psi_{\beta}e^{\langle \alpha,\gamma-\alpha\rangle\pi i}
\psi_{\alpha-\gamma}
\bar{\pi}_{\alpha}\psi_{\gamma-\alpha}(u)\nonumber\\
&=&e^{\langle \alpha,\gamma-\alpha\rangle \pi i}
\psi_{\alpha+\beta-\gamma}
\bar{\pi}_{\alpha}\psi_{\gamma-\alpha}(u).
\end{eqnarray}
On the other hand, we have:
\begin{eqnarray}
\bar{\pi}_{\alpha}\psi_{\beta}(u)
&=&e^{\langle \alpha,\gamma-\beta-\alpha\rangle\pi i}
\psi_{\alpha+\beta-\gamma}
\bar{\pi}_{\alpha}\psi_{\gamma-\beta-\alpha}\psi_{\beta}(u)\nonumber\\
&=&e^{\langle \alpha,\gamma-\beta-\alpha\rangle\pi i}
\psi_{\alpha+\beta-\gamma}\bar{\pi}_{\alpha}\psi_{\gamma-\alpha}(u).
\end{eqnarray}
The result follows. $\;\;\;\;\Box$

\begin{lem}\label{l3.10} {\it For any $\alpha\in L_{0}$
 we have}
\begin{eqnarray}
\bar{\pi}_{\alpha}(Y(u,z)v)=Y(u,z)\bar{\pi}_{\alpha}(v)\;\;\;
\mbox{{\it for any} }u, v\in U.
\end{eqnarray}
\end{lem}

{\bf Proof.} Without losing generality we assume that
$u\in V^{(\beta)}, v\in V^{(\gamma)}$, where
$\beta,\gamma\in L$. Then by Definition 3.3, (3.14) and Lemma 3.9 we obtain
\begin{eqnarray}
& &\;\;\;\;\bar{\pi}_{\alpha}(Y(u,z)v)\nonumber\\
& &=\bar{\pi}_{\alpha}\psi_{-\beta-\gamma}E^{-}(\beta,z)
Y(\psi_{\beta}\Delta(\gamma,z)(u),z)
\Delta(\beta,-z)\psi_{\gamma}(v)\nonumber\\
& &=e^{\langle \alpha,\beta\rangle\pi i}
\psi_{-\beta-\gamma}
E^{-}(\beta,z)Y(\psi_{\beta}\Delta(\gamma,z)(u),z)
\Delta(\beta,-z)\psi_{\gamma}\bar{\pi}_{\alpha}(v)\nonumber\\
& &=e^{\langle \alpha,\beta\rangle\pi i}
\psi_{\alpha-\beta-\gamma}\psi_{-\alpha}E^{-}(\beta,z)
Y(\psi_{\beta}\Delta(\gamma,z)(u),z)
\Delta(\beta,-z)\psi_{\gamma}\bar{\pi}_{\alpha}(v)\nonumber\\
& &=e^{\langle \alpha,\beta\rangle\pi i}
\psi_{\alpha-\beta-\gamma}E^{-}(\beta,z)\nonumber\\
& &\;\;\;\;\cdot Y(\Delta(-\alpha,z)\psi_{\beta}\Delta(\gamma,z)(u),z)
\psi_{-\alpha}\Delta(\beta,-z)\psi_{\gamma}\bar{\pi}_{\alpha}(v)\nonumber\\
& &=e^{\langle \alpha,\beta\rangle\pi i}
\psi_{\alpha-\beta-\gamma}E^{-}(\beta,z)z^{\langle \alpha,\beta\rangle}
(-z)^{-\langle \alpha,\beta\rangle}\nonumber\\
& &\;\;\;\;\cdot Y(\psi_{\beta}\Delta(\gamma-\alpha,z)(u),z)
\Delta(\beta,-z)\psi_{\gamma-\alpha}\bar{\pi}_{\alpha}(v)\nonumber\\
& &=\psi_{\alpha-\beta-\gamma}E^{-}(\beta,z)
Y(\psi_{\beta}\Delta(\gamma-\alpha,z)(u),z)
\Delta(\beta,-z)\psi_{\gamma-\alpha}\bar{\pi}_{\alpha}(v)\nonumber\\
& &=Y(u,z)\bar{\pi}_{\alpha}(v).
\end{eqnarray}
This proves the lemma. $\;\;\;\;\Box$

\begin{lem}\label{l3.11} {\it For any $u\in U, v\in V^{(\gamma)},
\alpha\in L_{0}, \gamma\in L$, we have}
\begin{eqnarray}
Y(\bar{\pi} _{\alpha}(u),z)v=e^{-\langle \alpha,\gamma\rangle\pi i}
\bar{\pi}_{\alpha}Y(u,z)v.
\end{eqnarray}
\end{lem}

{\bf Proof.}  By linearity we may assume that $u\in V^{(\beta)}$ for some
$\beta\in L$.
Using skew-symmetry, (3.15), Definition 3.3 and Lemma 3.9
we obtain
\begin{eqnarray}
& &\;\;\;\;Y(\bar{\pi} _{\alpha}(u),z)v\nonumber\\
& &=\psi_{\alpha-\beta-\gamma}E^{-}(\beta-\alpha,z)Y(\psi_{\beta-\alpha}
\Delta(\gamma,z)\bar{\pi}_{\alpha}(u),z)\Delta(\beta-\alpha,-z)\psi_{\gamma}(v)
\nonumber\\
& &=\psi_{\alpha-\beta-\gamma}E^{-}(\beta-\alpha,z)e^{zL(-1)}
Y(\Delta(\beta-\alpha,-z)\psi_{\gamma}(v),-z)\psi_{\beta-\alpha}
\Delta(\gamma,z)\bar{\pi}_{\alpha}(u)\nonumber\\
& &=\psi_{\alpha-\beta-\gamma}E^{-}(\beta-\alpha,z)e^{zL(-1)}\psi_{-\alpha}
Y(\Delta(\beta,-z)\psi_{\gamma}(v),-z)\psi_{\beta}
\Delta(\gamma,z)\bar{\pi}_{\alpha}(u)\nonumber\\
& &=e^{\langle \alpha,\beta\rangle\pi i}
\psi_{\alpha-\beta-\gamma}
E^{-}(\beta-\alpha,z)e^{zL(-1)}\psi_{-\alpha}\bar{\pi}_{\alpha}
Y(\Delta(\beta,-z)\psi_{\gamma}(v),-z)\psi_{\beta}
\Delta(\gamma,z)(u)\nonumber\\
& &=e^{\langle \alpha,\beta\rangle\pi i}
\psi_{\alpha-\beta-\gamma}
E^{-}(\beta-\alpha,z)e^{zL(-1)}\psi_{-\alpha}\bar{\pi}_{\alpha}\nonumber\\
& &\;\;\;\;\cdot
e^{-zL(-1)}Y(\psi_{\beta}\Delta(\gamma,z)(u),z)\Delta(\beta,-z)
\psi_{\gamma}(v)\nonumber\\
& &=e^{\langle \alpha,\beta\rangle\pi i}
\psi_{\alpha-\beta-\gamma}
E^{-}(\beta-\alpha,z)e^{zL(-1)}\psi_{-\alpha}e^{-zL(-1)}\bar{\pi}_{\alpha}\nonumber\\
& &\;\;\;\;\cdot Y(\psi_{\beta}\Delta(\gamma,z)(u),z)\Delta(\beta,-z)
\psi_{\gamma}(v)\nonumber\\
& &=e^{\langle \alpha,\beta\rangle\pi i}
\psi_{\alpha-\beta-\gamma}
E^{-}(\beta-\alpha,z)E^{-}(\alpha,z)\psi_{-\alpha}\bar{\pi}_{\alpha}
Y(\psi_{\beta}\Delta(\gamma,z)(u),z)\Delta(\beta,-z)
\psi_{\gamma}(v)\nonumber\\
& &=e^{\langle \alpha,\beta\rangle\pi i}
\psi_{\alpha-\beta-\gamma}
E^{-}(\beta,z)\psi_{-\alpha}\bar{\pi}_{\alpha}Y(\psi_{\beta}\Delta(\gamma,z)(u),z)
\Delta(\beta,-z)\psi_{\gamma}(v)\nonumber\\
& &=e^{\langle \alpha,\beta\rangle\pi i}\psi_{-\beta-\gamma}
E^{-}(\beta,z)\bar{\pi}_{\alpha}Y(\psi_{\beta}\Delta(\gamma,z)(u),z)
\Delta(\beta,-z)\psi_{\gamma}(v)\nonumber\\
& &=e^{\langle \alpha,\beta\rangle\pi i}
e^{-\langle \alpha,\beta+\gamma\rangle\pi i}
\bar{\pi}_{\alpha}\psi_{-\beta-\gamma}
E^{-}(\beta,z)Y(\psi_{\beta}\Delta(\gamma,z)(u),z)
\Delta(\beta,-z)\psi_{\gamma}(v)\nonumber\\
& &=e^{-\langle \alpha,\gamma\rangle\pi i}\bar{\pi}_{\alpha}Y(u,z)v.\ \ \Box
\end{eqnarray}

\begin{rem}\label{r3.12}
Let $I$ be the sum of all subspaces $(\bar{\pi}_{\alpha}-1)U$ of $U$ for
$\alpha\in L_{0}$. Define $\bar{U}$ to be the quotient space $U/I.$
Then the multiplicity of the $\sigma_{\beta}$-twisted $V$-module $V^{(\beta)}$
in $\bar{U}$ is exactly one for $\beta\in L$.
It follows from Lemma 3.10 that $I$ is
 a left ideal of the generalized vertex  algebra $U$.
But from Lemma 3.11, $I$ is not necessarily
a right ideal, {\it i.e.,} for $u\in I, v\in U$, $Y(u,z)v$ may not be in
$I\{z\}$
unless $\langle \alpha,\beta\rangle\in 2{\Z}$ for any $\alpha\in L_{0},
\beta\in L$ (Lemma 3.11). In general, $\bar{U}$ is a $U$-module, but it is not
a
quotient generalized vertex algebra of $U$.
\end{rem}

Before we modify the definition of vertex operator $Y(\cdot,z)$ (3.23)
to get an
abelian intertwining algebra we consider a special case.
Let $L_{1}$ be an integral sublattice of $L$ such that $L_{0}\subseteq
L_{1}\subseteq L$ and
\begin{eqnarray}
\langle \lambda,\beta\rangle \in {\Z},\;
\langle \alpha,\beta\rangle \in 2{\Z}\;\;\;\mbox{ for any }\lambda\in
P,\alpha\in L_{0},
\beta\in L_{1}.\label{3.56}
\end{eqnarray}
Set $U_{1}=\oplus_{\beta\in L_{1}}V^{(\beta)}$ and $\bar{U}_{1}=U_{1}/I$.
Then it follows from Theorem 3.5, Lemmas 3.10 and 3.11 that $\bar{U}_{1}$ is a
generalized vertex algebra. By Lemma 3.8, for any $\alpha_{j}\in
L_{1},\lambda_{j}\in P$
for $j=1,2$ we have:
\begin{eqnarray}
& &\eta((\alpha_{1},\lambda_{1}),(\alpha_{2},\lambda_{2}))
=-\langle \alpha_{1},\alpha_{2}\rangle-\langle \alpha_{1},\lambda_{2}\rangle
-\langle \alpha_{2},\lambda_{1}\rangle\in {\Z}/2{\Z},\\
& &C((\alpha_{1},\lambda_{1}),(\alpha_{2},\lambda_{2}))
=(-1)^{\langle \alpha_{1},\lambda_{2}\rangle
+\langle \alpha_{2},\lambda_{1}\rangle}.
\end{eqnarray}
Then
\begin{eqnarray}
& &\;\;\;C((\alpha_{1},\lambda_{1}),(\alpha_{2},\lambda_{2}))
z_{0}^{-1}\delta\left(\frac{z_{2}-z_{1}}{-z_{0}}\right)
\left(\frac{z_{2}-z_{1}}{z_{0}}\right)^{\eta((\alpha_{1},\lambda_{1}),(\alpha_{2},\lambda_{2}))}
\nonumber\\
& &=(-1)^{\langle \alpha_{1},\alpha_{2}\rangle}
z_{0}^{-1}\delta\left(\frac{z_{2}-z_{1}}{-z_{0}}\right).
\end{eqnarray}
Therefore for any $u\in V^{(\alpha_{1},\lambda_{1})},
v\in V^{(\alpha_{2},\lambda_{2})}, w\in V^{(\alpha_{3},\lambda_{3})},
(\alpha_{j},\alpha_{j})\in L_{1}\times P$, the generalized Jacobi identity
(3.26) becomes
 the following super Jacobi identity:
\begin{eqnarray}
& &z_{0}^{-1}\delta\left(\frac{z_{1}-z_{2}}{z_{0}}\right)
Y(u,z_{1})Y(v,z_{2})w\nonumber\\
&-&(-1)^{\langle
\alpha_{1},\alpha_{2}\rangle}z_{0}^{-1}\delta\left(\frac{z_{2}-z_{1}}
{-z_{0}}\right)
Y(v,z_{2})Y(u,z_{1})w\nonumber\\
& &=z_{2}^{-1}\delta\left(\frac{z_{1}-z_{0}}{z_{2}}\right)
Y(Y(u,z_{0})v,z_{2})w.
\end{eqnarray}
Then we have:

\bc{c3.13}  {\it Let $L_{1}$ be an integral sublattice of $L$
satisfying (\ref{3.56}).
Then $U_{1}$ is a vertex superalgebra with $I$ as an ideal so that
$\bar{U}_{1}$ is
a quotient vertex superalgebra.}
\ec

Continuing with Corollary 3.13, let $M$ be an irreducible $V$-module
such that $\alpha(0)$ has rational eigenvalues on $M$ for any $\alpha\in
L_{1}$.
Let $\gamma$ be an $H$-weight of $M$. Then $P(M)=\gamma+P$.
Set $W_{1}=\oplus _{\alpha\in L_{1}}M^{(\alpha)}$.
For any $(\alpha_{1},\lambda_{1})\in L_{1}\times P,
(\alpha_{3},\lambda_{3})\in L_{1}\times P(M)=L_{1}\times (\gamma+P)$,
since $\langle \alpha_{1},\alpha_{3}\rangle, \langle
\alpha_{3},\lambda_{1}\rangle
\in {\Z}$ and $\langle \alpha_{1},\lambda\rangle\in {\Z}$ for any $\lambda\in
P$,
we have:
\begin{eqnarray}
& &\;\;\;\;\eta((\alpha_{1},\lambda_{1}),(\alpha_{3},\lambda_{3}))\nonumber\\
& &=-\langle \alpha_{1},\alpha_{3}\rangle -\langle
\alpha_{1},\lambda_{3}\rangle
-\langle \alpha_{3},\lambda_{1}\rangle\nonumber\\
& &=-\langle \alpha_{1},\gamma\rangle\in {1\over T}{\Z}/{\Z}.
\end{eqnarray}
Then we have the following twisted Jacobi identity:
\begin{eqnarray}
& &z_{0}^{-1}\delta\left(\frac{z_{1}-z_{2}}{z_{0}}\right)
Y(u,z_{1})Y(v,z_{2})w\nonumber\\
&-&(-1)^{\langle
\alpha_{1},\alpha_{2}\rangle}z_{0}^{-1}\delta\left(\frac{z_{2}-z_{1}}
{-z_{0}}\right)
Y(v,z_{2})Y(u,z_{1})w\nonumber\\
& &=z_{2}^{-1}\delta\left(\frac{z_{1}-z_{0}}{z_{2}}\right)
\left(\frac{z_{2}+z_{0}}{z_{1}}\right)^{-\langle \gamma,\alpha_{1}\rangle}
Y(Y(u,z_{0})v,z_{2})w.
\end{eqnarray}
Therefore, for any $u\in V^{(\alpha_{1},\lambda_{1})},
v\in V^{(\alpha_{2},\lambda_{2})}, w\in M^{(\alpha_{3},\lambda_{3})},
(\alpha_{1},\lambda_{1}),(\alpha_{2},\lambda_{2})\in L_{1}\times P,
(\alpha_{3},\lambda_{3})
\in L_{1}\times P(M)$, we have the following
super Jacobi identity:
\begin{eqnarray}
& &z_{0}^{-1}\delta\left(\frac{z_{1}-z_{2}}{z_{0}}\right)
Y(u,z_{1})Y(v,z_{2})w\nonumber\\
&-&(-1)^{\langle
\alpha_{1},\alpha_{2}\rangle}z_{0}^{-1}\delta\left(\frac{z_{2}-z_{1}}
{-z_{0}}\right)
Y(v,z_{2})Y(u,z_{1})w\nonumber\\
& &=z_{2}^{-1}\delta\left(\frac{z_{1}-z_{0}}{z_{2}}\right)
\left(\frac{z_{2}+z_{0}}{z_{1}}\right)^{\eta((\alpha_{1},\lambda_{1}),(\alpha_{3},\lambda_{3}))}
Y(Y(u,z_{0})v,z_{2})w.
\end{eqnarray}
It is clear that $\sigma_{\gamma}=e^{-2\pi i \gamma(0)}$ is an
automorphism of $\bar{U}_{1}$. Then $W_{1}$ is a $\sigma_{\gamma}$-twisted
$\bar{U}_{1}$-module with $I(W_{1})$ as a submodule, where $I(W_{1})$ is
defined as the subspace
linearly spanned by $(\pi_{\alpha}-1)W_{1}$ for $\alpha\in L_{0}$.
Summarizing the previous arguments we have:

\bc{c3.14}  {\it Let $L_{1}$ be an integral sublattice of $L$
satisfying (\ref{3.56}). Let $M$ be an irreducible $V$-module
such that $\alpha(0)$ has rational eigenvalues on $M$ for any $\alpha\in
L_{1}$.
Then $W_{1}=\oplus_{\beta\in L_{1}}M^{(\beta)}$ with $M^{(0)}=M$, defined as in
Theorem 3.6,
is a $\sigma_{\gamma}$-twisted $\bar{U}_{1}$-module with a submodule
$I(W_{1})$, so that
$\bar{W}_{1}=W_{1}/I(W_{1})$ is a quotient $\sigma_{\gamma}$-twisted
$\bar{U}_{1}$-module.}
\ec

\bc{c3.15} {\it Under the conditions of Corollary 3.14, assume that there is a
$\tau$-twisted
$\bar{U}_{1}$-module $E$ containing $M$ as a $V$-submodule for some
finite-order automorphism
$\tau$ of $\bar{U}_{1}$. Then $\tau=\sigma_{\gamma}$.}
\ec

{\bf Proof.} Since $\bar{U}_{1}$ is simple and $V^{(\alpha)}$ is a simple
current for
any $\alpha\in L_{1}$,
$(M^{(\alpha)},Y_{\bar{W}_{1}}(\cdot,z))$ is a tensor product of $V^{(\alpha)}$
with $M$. Let
$E^{\alpha}$ be the subspace of $E$ linearly spanned by $u_{n}M$ for $u\in
V^{(\alpha)}, n\in {\Q}$.
Similarly, $(E^{\alpha},Y_{E}(\cdot,z))$ is also a tensor product of
$V^{(\alpha)}$ with $M$.
Therefore, there is a $c_{\alpha}\in {\C}^{*}$ such that
$Y_{\bar{W}_{1}}(u,z)w=c_{\alpha}Y_{E}(u,z)w$
for any $u\in V^{(\alpha)}, w\in M$. By definition of a twisted module,
$\tau$ and $\sigma_{\gamma}$ have the same order.
Then $\tau|_{V^{(\alpha)}}=\sigma_{\gamma}|_{V^{(\alpha)}}$
for any $\alpha\in L_{1}$. Thus $\tau=\sigma_{\gamma}$.\ \  $\Box$

As mentioned in Remark 3.12, in general $\bar{U}=U/I$ is not a generalized
vertex algebra. Next
we shall modify the definition (3.23) and prove that $\bar{U}=U/I$ is an
abelian intertwining
algebra.

For any $\alpha,\beta\in L_{0}$, we have a $V$-isomorphism
$\bar{\pi}_{\alpha}\bar{\pi}_{\beta}\bar{\pi}_{\alpha+\beta}^{-1}$ from $V$
onto $V$. Since $V$ is a simple vertex operator algebra,
by Schur's lemma there is a nonzero complex number $A_{0}(\alpha,\beta)$ such
that
\be{3.62}
\bar{\pi}_{\alpha+\beta}
=A_{0}(\alpha,\beta)\bar{\pi}_{\alpha}\bar{\pi}_{\beta},
\ee
where both sides are considered as $V$-isomorphisms from
$V^{(\alpha+\beta)}$ onto $V$. For any $\gamma\in L$ and for any
$u\in V^{(\gamma)}$, we have:
\begin{eqnarray}
\bar{\pi}_{\alpha+\beta}(u)
&=&\bar{\pi}_{\alpha+\beta}\psi_{\alpha+\beta-\gamma}
\psi_{\gamma-\alpha-\beta}(u)\nonumber\\
&=&e^{-\langle \alpha+\beta,\alpha+\beta-\gamma\rangle\pi i}
\psi_{\alpha+\beta-\gamma}\bar{\pi}_{\alpha+\beta}\psi_{\gamma-\alpha-\beta}(u)
\nonumber\\
&=&e^{-\langle \alpha+\beta,\alpha+\beta-\gamma\rangle\pi i}
A_{0}(\alpha,\beta)\psi_{\alpha+\beta-\gamma}\bar{\pi}_{\alpha}
\bar{\pi}_{\beta}
\psi_{\gamma-\alpha-\beta}(u)\nonumber\\
&=&A_{0}(\alpha,\beta)\bar{\pi}_{\alpha}\bar{\pi}_{\beta}(u).
\end{eqnarray}
Then (\ref{3.62})  holds when both sides are considered as operators on $U$.
It is easy to see that the following $2$-cocycle condition hold:
\begin{eqnarray}
A_{0}(\alpha_{1}+\alpha_{2},\alpha_{3})A_{0}(\alpha_{1},\alpha_{2})
=A_{0}(\alpha_{1},\alpha_{2}+\alpha_{3})A_{0}(\alpha_{2},\alpha_{3})
\end{eqnarray}
for any $\alpha_{i}\in L_{0}, i=1,2,3$.

Next we define
\begin{eqnarray}
C_{0}(\alpha,\beta)=A_{0}(\alpha,\beta)A_{0}(\beta,\alpha)^{-1}\;\;\;
\mbox{for }\alpha,\beta\in L_{0}.
\end{eqnarray}
Then  $C_{0}(\cdot,\cdot)$ satisfies the  properties (\ref{3.22}) and
\begin{eqnarray}
\bar{\pi}_{\beta}\bar{\pi}_{\alpha}
=C_{0}(\alpha,\beta)\bar{\pi}_{\alpha}\bar{\pi}_{\beta}
\;\;\;\mbox{ for }\alpha,\beta\in L_{0}.
\end{eqnarray}

Since $L_{0}$ is a sublattice of $L$, there is a basis
$\{\beta_{1},\beta_{2},\cdots,\beta_{n}\}$ for $L$ and a basis
$\{\alpha_{1},\alpha_{2},\cdots,\alpha_{n}\}$ for
$L_{0}$ such that each $\alpha_{i}$ is an integral multiple of $\beta_{i}$.
It is easy to find a
${\Z}$-bilinear function $A_{1}(\cdot,\cdot)$ on $L$ with values in
${\C}^{*}$ satisfying the following condition:
\begin{eqnarray}
A_{0}(\alpha_{i},\alpha_{j})=A_{1}(\alpha_{i},\alpha_{j})^{2}
\;\;\;\mbox{for any }1\le i,j\le n.
\end{eqnarray}
Fixing such an $A_{1}(\cdot,\cdot)$, we define $C_{1}(\cdot,\cdot)$ on $L\times
L$ as
follows:
\begin{eqnarray}
C_{1}(\alpha,\beta)=A_{1}(\alpha,\beta)A_{1}(\beta,\alpha)^{-1}\;\;\;\mbox{for
any }
\alpha,\beta\in L.
\end{eqnarray}
Then $C_{1}(\cdot,\cdot)$ satisfies the following conditions:
\begin{eqnarray}
C_{1}(\beta,\beta)=1,\;C_{1}(\beta_{1},\beta_{2})C_{1}(\beta_{2},\beta_{1})=1,\;
C_{1}(\beta_{1}+\beta_{2},\beta_{3})=C_{1}(\beta_{1},\beta_{3})C_{1}(\beta_{2},\beta_{3})
\end{eqnarray}
for any $\beta,\beta_{1},\beta_{2},\beta_{3}\in L$.

Next, for any $\alpha\in L_{0}$, we define a linear automorphism
$\pi_{\alpha}$ on $U$ as follows:
\begin{eqnarray}
\pi_{\alpha}(u)=C_{1}(\beta,\alpha)\bar{\pi}_{\alpha}(u)\;\;\;\mbox{for any }
u\in V^{(\beta)}\subseteq U.
\end{eqnarray}
Then for any $\alpha_{1},\alpha_{2}\in L_{0}, \beta\in L$, we have:
\begin{eqnarray}
\pi_{\alpha_{1}}\pi_{\alpha_{2}}(u)
&=&C_{1}(\beta-\alpha_{2},\alpha_{1})C_{1}(\beta,\alpha_{2})
\bar{\pi}_{\alpha_{1}}\bar{\pi}_{\alpha_{2}}(u)\nonumber\\
&=&C_{1}(\beta-\alpha_{2},\alpha_{1})C_{1}(\beta,\alpha_{2})
C_{0}(\alpha_{2},\alpha_{1})\bar{\pi}_{\alpha_{2}}\bar{\pi}_{\alpha_{1}}(u)
\nonumber\\
&=&C_{1}(\beta-\alpha_{2},\alpha_{1})C_{1}(\beta,\alpha_{2})
C_{0}(\alpha_{2},\alpha_{1})C_{1}(\alpha_{2},\beta-\alpha_{1})
C_{1}(\alpha_{1},\beta)
\pi_{\alpha_{2}}\pi_{\alpha_{1}}(u)\nonumber\\
&=&\pi_{\alpha_{2}}\pi_{\alpha_{1}}(u)
\end{eqnarray}
for any $u\in V^{(\beta)}$. Thus
\begin{eqnarray}
\pi_{\alpha_{1}}\pi_{\alpha_{2}}=\pi_{\alpha_{2}}\pi_{\alpha_{1}}
\;\;\;\mbox{for any }\alpha_{1},\alpha_{2}\in L_{0}.
\end{eqnarray}

\begin{rem}\label{r3.16} If $L={\Z}\alpha$ is of rank one, then $L_{0}=kL$ for
some
positive integer $k$. We can fix $\bar{\pi}_{k\alpha}$ first, then we define
$\bar{\pi}_{nk\alpha}=\bar{\pi}_{k\alpha}^{n}$ for any $n\in {\Z}$. Then
$C_{0}(\cdot,\cdot)\equiv 1$. So we can take $C_{1}(\cdot,\cdot)\equiv 1$.
\end{rem}

\begin{lem}\label{l3.17} {\it For any $\alpha\in L_{0}, \beta\in L$, we have}
\begin{eqnarray}
\psi_{\beta}\pi_{\alpha}
=e^{\langle \alpha,\beta\rangle \pi
i}C_{1}(\beta,\alpha)\pi_{\alpha}\psi_{\beta}.
\end{eqnarray}
\end{lem}

{\bf Proof.} Let $\gamma\in L, u\in V^{(\gamma)}$. Then by the definition
(3.71) of
 $\pi_{\alpha}$ we have:
\begin{eqnarray}
\psi_{\beta}\pi_{\alpha}(u)
&=&C_{1}(\gamma,\alpha)\psi_{\beta}\bar{\pi}_{\alpha}(u)\nonumber\\
&=&C_{1}(\gamma,\alpha)e^{\langle \alpha,\beta\rangle\pi i}
\bar{\pi}_{\alpha}\psi_{\beta}(u)\nonumber\\
&=&C_{1}(\gamma,\alpha)C_{1}(\alpha,\gamma-\beta)e^{\langle
\alpha,\beta\rangle\pi i}
\pi_{\alpha}\psi_{\beta}(u)\nonumber\\
&=&e^{\langle \alpha,\beta\rangle\pi
i}C_{1}(\beta,\alpha)\pi_{\alpha}\psi_{\beta}(u).
\;\;\;\;\Box
\end{eqnarray}

\begin{lem}\label{l3.18} {\it For any $\alpha\in L_{0}$
 we have}
\begin{eqnarray}
\pi_{\alpha}(Y(u,z)v)=C_{1}(\beta,\alpha)Y(u,z)\pi_{\alpha}(v)\;\;\;
\mbox{{\it for any} }u\in V^{(\beta)}\subseteq U.
\end{eqnarray}
\end{lem}

{\bf Proof.} Let $u\in V^{(\beta)}, v\in V^{(\gamma)}$ with
$\beta,\gamma\in L$. Then by the definition (3.71) we have:
\begin{eqnarray}
\pi_{\alpha}(Y(u,z)v)
&=&C_{1}(\beta+\gamma,\alpha)\bar{\pi}_{\alpha}(Y(u,z)v)\nonumber\\
&=&C_{1}(\beta+\gamma,\alpha)Y(u,z)\bar{\pi}_{\alpha}v\nonumber\\
&=&C_{1}(\beta+\gamma,\alpha)C_{1}(\gamma,\alpha)^{-1}Y(u,z)\pi_{\alpha}v\nonumber\\
&=&C_{1}(\beta,\alpha)Y(u,z)\pi_{\alpha}(v).
\end{eqnarray}
This proves the assertion. $\;\;\;\;\Box$

\begin{lem}\label{l3.19} {\it For any $u\in V^{(\beta)}, v\in V^{(\gamma)},
\alpha\in L_{0}, \gamma\in L$, we have}
\begin{eqnarray}
Y(\pi _{\alpha}(u),z)v=e^{-\langle \alpha,\gamma\rangle\pi i}
C(\alpha,\gamma)\pi_{\alpha}Y(u,z)v.
\end{eqnarray}
\end{lem}

{\bf Proof.} We may assume that $u\in V^{(\beta)}$ for some $\beta\in L$. Then
by definition we have:
\begin{eqnarray}
Y(\pi _{\alpha}(u),z)v
&=&C_{1}(\beta,\alpha)Y(\bar{\pi}_{\alpha}(u),z)v\nonumber\\
&=&e^{-\langle \alpha,\gamma\rangle\pi i}C_{1}(\beta,\alpha)
\bar{\pi}_{\alpha}Y(u,z)v\nonumber\\
&=&e^{-\langle \alpha,\gamma\rangle\pi i}C_{1}(\beta,\alpha)
C_{1}(\beta+\gamma,\alpha)^{-1}\pi_{\alpha}Y(u,z)v\nonumber\\
&=&e^{-\langle \alpha,\gamma\rangle\pi i}C_{1}(\alpha,\gamma)
\pi_{\alpha}Y(u,z)v.
\end{eqnarray}
Then the proof is complete. $\;\;\;\;\Box$

\begin{rem}\label{r3.20} Let $\alpha_{1},\cdots,\alpha_{n}$ be a basis of
$L_{0}$.
Then we can define a ${\Z}$-linear map $\pi'$ from $L_{0}$ to $Aut (U)$ as
follows:
\begin{eqnarray}
\pi'_{\alpha}=\pi_{\alpha_{1}}^{k_{1}}\pi_{\alpha_{2}}^{k_{2}}\cdots
\pi_{\alpha_{n}}^{k_{n}}
\end{eqnarray}
for any
$\alpha=k_{1}\alpha_{1}+k_{2}\alpha_{2}+\cdots +k_{n}\alpha_{n}\in L_{0}$.
Since all $\pi_{\alpha_{i}}$'s commute each other, $\pi'$ is well-defined. It
is easy to see that
\begin{eqnarray}
\pi'_{\alpha+\beta}=\pi'_{\alpha}\pi'_{\beta}=\pi'_{\beta}\pi'_{\alpha}
\;\;\;\mbox{ for any }\alpha,\beta\in L_{0}.
\end{eqnarray}
It is also easy to see that Lemmas 3.17, 3.18 and 3.19 still hold. By slightly
abusing the notion,  from now on we will use $\pi'$ for $\pi$.
\end{rem}

Let $\bar{A}=L\times P/D$, where $D=\{ (\alpha,-\alpha)|\alpha\in L_{0}\}$
is a subgroup of $A=L\times P$. Let $\{ \lambda_{i}|i\in \bar{A}\}$ be a
(complete) set of
representatives in $A.$
Then we define $\bar{V}=\oplus_{i\in \bar{A}}V^{(\lambda_{i})}$.
For any $u\in V^{(\lambda_{i})},v\in V^{(\lambda_{j})}$, we define
$\bar{Y}(u,z)v\in V^{(\lambda_{i+j})}\{z\}$ as follows:
\begin{eqnarray}
\bar{Y}(u,z)v
=C_{1}(\lambda_{j},\lambda_{i})\pi_{\lambda_{i}+\lambda_{j}-\lambda_{i+j}}Y(u,z)v.
\end{eqnarray}
Because $L(-1)$ commutes with $\pi_{\alpha}$ for $\alpha\in L_{0}$ and
$L(-1)V^{(\beta)}
\subseteq V^{(beta)}$ for $\beta\in L$, the $L(-1)$-derivative property
(Proposition 3.4)
still holds.

We define a function $h(i,j,k)$ from $\bar{A}\times \bar{A}\times \bar{A}$ to
${\C}^{*}$
as follows:
\begin{eqnarray}
h(i,j,k)
=e^{-\langle \lambda_{i}+\lambda_{j}-\lambda_{i+j},\lambda_{k}\rangle \pi i}
C_{1}(\lambda_{i}+\lambda_{j}-\lambda_{i+j},\lambda_{k})^{2}\;\;\;\mbox{for
}i,j,k\in \bar{A}.
\end{eqnarray}
Next we shall prove that  $h(i,j,k)$ is a $3$-cocycle. Observe that
$h(i,j,k)$ is symmetric in the first two variables with
the third variable being fixed. Following [DL], we prove that $h(i,j,k)$ is a
$3$-cocycle
by proving that with the third variable $k$ being fixed,
$h(i,j,k)$ is a $2$-cocycle, {\it i.e.,}
\begin{eqnarray}
h(i,j,k)h(i,j+r,k)^{-1}h(i+j,r,k)h(j,r,k)^{-1}=1\;\;\;\mbox{for }i,j,r\in
\bar{A}.
\end{eqnarray}
For $i,j,r,k\in \bar{A}$ we have:
\begin{equation}\label{3.86}
e^{-\langle \lambda_{i}+\lambda_{j}-\lambda_{i+j},\lambda_{k}\rangle \pi i}
e^{\langle \lambda_{i}+\lambda_{j+r}-\lambda_{i+j+r},\lambda_{k}\rangle \pi i}
=e^{\langle \lambda_{i+j}+\lambda_{r}-\lambda_{i+j+r},\lambda_{k}\rangle \pi i}
e^{\langle \lambda_{j}+\lambda_{r}-\lambda_{j+r},\lambda_{k}\rangle \pi i}
\end{equation}
and
\begin{eqnarray}
& &\ \ C_{0}(\lambda_{i}+\lambda_{j}-\lambda_{i+j},\lambda_{k})
C_{0}(\lambda_{i}+\lambda_{j+r}-\lambda_{i+j+r},\lambda_{k})^{-1}\nonumber\\
& &= C_{0}(\lambda_{i+j}+\lambda_{r}-\lambda_{i+j+r},\lambda_{k})
C_{0}(\lambda_{j}+\lambda_{r}-\lambda_{j+r},\lambda_{k}).\label{3.87}
\end{eqnarray}
Combining (3.85) with (3.86) we obtain (3.84).
Then from Proposition 12.13 of
[DL], $h(i,j,k)$ is a $3$-cocycle. Next we define
\begin{eqnarray}
\bar{C}((\lambda_{1},h_{1}),(\lambda_{2},h_{2}))
=C((\lambda_{1},h_{1}),(\lambda_{2},h_{2}))C_{1}(\lambda_{1},\lambda_{2})
\end{eqnarray}
for any $(\lambda_{i}\times P)\in L\times P$ (recall the definition of
$C(\cdot,\cdot)$ from (3.20)).
 Then $\bar{C}$ satisfies (3.22).

\bt{t3.21} {\it For any $u\in V^{(\lambda_{i},h_{1})}, v\in
V^{(\lambda_{j},h_{2})},
w\in V^{(\lambda_{k},h_{3})}$, the following generalized Jacobi identity
holds:}
\begin{eqnarray}
& &z_{0}^{-1}\delta\left(\frac{z_{1}-z_{2}}{z_{0}}\right)
\left(\frac{z_{1}-z_{2}}{z_{0}}\right)^{\eta((\lambda_{i},h_{1}),(\lambda_{j},h_{2}))}
\bar{Y}(u,z_{1})\bar{Y}(v,z_{2})w\nonumber\\
&-&\bar{C}((\lambda_{i},h_{1}),(\lambda_{j},h_{2}))
z_{0}^{-1}\delta\left(\frac{z_{2}-z_{1}}{-z_{0}}\right)
\left(\frac{z_{2}-z_{1}}{z_{0}}\right)^{\eta((\lambda_{i},h_{1}),(\lambda_{j},h_{2}))}
\nonumber\\
& &\cdot\bar{Y}(v,z_{2})\bar{Y}(u,z_{1})w\nonumber\\
&=&z_{2}^{-1}\delta\left(\frac{z_{1}-z_{0}}{z_{2}}\right)
\left(\frac{z_{2}+z_{0}}{z_{1}}\right)^{\eta((\lambda_{i},h_{1}),(\lambda_{k},h_{3}))}
h(i,j,k)\bar{Y}(\bar{Y}(u,z_{0})v,z_{2})w.
\end{eqnarray}
{\it Therefore $(\bar{U}, {\bf 1}, \omega, \bar{Y}, T, \bar{A}, \eta, \bar{C})$
is an
abelian intertwining algebra.}
\et

{\bf Proof.} Since $Y(v,z)w\in V^{(\lambda_{i}+\lambda_{j},
h_{2}+h_{3})}\{z\}$,
using (3.82) and Lemma 3.18 we obtain
\begin{eqnarray}
& &\bar{Y}(u,z_{1})\bar{Y}(v,z_{2})w\nonumber\\
&=&C_{1}(\lambda_{j+k},\lambda_{i})\pi_{\lambda_{i}+\lambda_{j+k}-\lambda_{i+j+k}}
Y(u,z_{1})\bar{Y}(v,z_{2})w,\nonumber\\
&=&C_{1}(\lambda_{j+k},\lambda_{i})C_{1}(\lambda_{k},\lambda_{j})
\pi_{\lambda_{i}+\lambda_{j+k}-\lambda_{i+j+k}}
Y(u,z_{1})\pi_{\lambda_{j}+\lambda_{k}
-\lambda_{j+k}}Y(v,z_{2})w\nonumber\\
&=&C_{1}(\lambda_{j+k},\lambda_{i})C_{1}(\lambda_{k},\lambda_{j})
C_{1}(\lambda_{j}+\lambda_{k}-\lambda_{j+k},\lambda_{i})
\pi_{\lambda_{i}+\lambda_{j+k}-\lambda_{i+j+k}}
\pi_{\lambda_{i}+\lambda_{j+k}-\lambda_{i+j+k}}\nonumber\\
& &\cdot \pi_{\lambda_{j}+\lambda_{k}
-\lambda_{j+k}}Y(u,z_{1})Y(v,z_{2})w\nonumber\\
&=&C_{1}(\lambda_{k},\lambda_{j})C_{1}(\lambda_{j}+\lambda_{k},\lambda_{i})
\pi_{\lambda_{i}+\lambda_{j}+\lambda_{k}-\lambda_{i+j+k}}
Y(u,z_{1})Y(v,z_{2})w.
\end{eqnarray}
Symmetrically, we have:
\begin{eqnarray}
& &\ \ \ \bar{Y}(v,z_{2})\bar{Y}(u,z_{1})w\nonumber\\
& &=C_{1}(\lambda_{k},\lambda_{i})C_{1}(\lambda_{i}+\lambda_{k},\lambda_{j})
\pi_{\lambda_{i}+\lambda_{j}+\lambda_{k}-\lambda_{i+j+k}}
Y(v,z_{2})Y(u,z_{1})w.
\end{eqnarray}
On the other hand, using (3.82) and  Lemma 3.19 we get:
\begin{eqnarray}
& &\bar{Y}(\bar{Y}(u,z_{0})v,z_{2})w\nonumber\\
&=&C_{1}(\lambda_{k},\lambda_{j+k})\pi_{\lambda_{i+j}+\lambda_{k}-\lambda_{i+j+k}}
Y(\bar{Y}(u,z_{0})v,z_{2})w\nonumber\\
&=&C_{1}(\lambda_{k},\lambda_{j+k})C_{1}(\lambda_{j},\lambda_{i})
\pi_{\lambda_{i+j}+\lambda_{k}-\lambda_{i+j+k}}
Y(\pi_{\lambda_{i}+\lambda_{j}-\lambda_{i+j}}Y(u,z_{0})v,z_{2})w\nonumber\\
&=&C_{1}(\lambda_{k},\lambda_{j+k})C_{1}(\lambda_{j},\lambda_{i})
C_{1}(\lambda_{i}+\lambda_{j}-\lambda_{i+j},\lambda_{k})
e^{-\langle \lambda_{i}+\lambda_{j}-\lambda_{i+j},\lambda_{k}\rangle \pi i}
\nonumber\\
& &\cdot \pi_{\lambda_{i+j}+\lambda_{k}-\lambda_{i+j+k}}
\pi_{\lambda_{i}+\lambda_{j}-\lambda_{i+j}}Y(Y(u,z_{0})v,z_{2})w\nonumber\\
&=&C_{1}(\lambda_{k},\lambda_{j+k})C_{1}(\lambda_{j},\lambda_{i})
C_{1}(\lambda_{i}+\lambda_{j}-\lambda_{i+j},\lambda_{k})
e^{-\langle \lambda_{i}+\lambda_{j}-\lambda_{i+j},\lambda_{k}\rangle \pi i}
\nonumber\\
& &\cdot \pi_{\lambda_{i}+\lambda_{j}+\lambda_{k}-\lambda_{i+j+k}}
Y(Y(u,z_{0})v,z_{2})w.
\end{eqnarray}
Multiplying the generalized Jacobi identity (\ref{3.31}) by
$C_{1}(\lambda_{k},\lambda_{j})
C_{1}(\lambda_{j}+\lambda_{k},\lambda_{i})$, applying
$\pi_{\lambda_{i}+\lambda_{j}+\lambda_{k}-\lambda_{i+j+k}}$, then using
(3.89)-(3.91)
 we obtain
\begin{eqnarray}
& &z_{0}^{-1}\delta\left(\frac{z_{1}-z_{2}}{z_{0}}\right)
\left(\frac{z_{1}-z_{2}}{z_{0}}\right)^{\eta((\alpha,h_{1}),(\beta,h_{2}))}
\bar{Y}(u,z_{1})\bar{Y}(v,z_{2})w\nonumber\\
&-&C((\lambda_{i},h_{1}),(\lambda_{j},h_{2}))C_{1}(\lambda_{i},\lambda_{j})^{-2}
z_{0}^{-1}\delta\left(\frac{z_{2}-z_{1}}
{-z_{0}}\right)
\left(\frac{z_{2}-z_{1}}{z_{0}}\right)^{\eta((\lambda_{i},h_{1}),(\lambda_{j},h_{2}))}
\nonumber\\
& &\cdot \bar{Y}(v,z_{2})\bar{Y}(u,z_{1})w\nonumber\\
&=&e^{-\langle \lambda_{i}+\lambda_{j}-\lambda_{i+j},\lambda_{k}\rangle \pi i}
C_{1}(\lambda_{i}+\lambda_{j}-\lambda_{i+j},\lambda_{k})^{2}\cdot\nonumber\\
& &\cdot
z_{2}^{-1}\delta\left(\frac{z_{1}-z_{0}}{z_{2}}\right)
\left(\frac{z_{2}+z_{0}}{z_{1}}\right)^{\eta((\lambda_{i},h_{1}),(\lambda_{k},h_{3}))}
\bar{Y}(\bar{Y}(u,z_{0})v,z_{2})w.
\end{eqnarray}
This gives the generalized Jacobi identity (3.88). The proof is complete.
$\;\;\;\;\Box$

\section{Rationality for certain extensions of vertex operator algebras}

In this section we study the rationality of certain  extensions of
rational vertex operator algebras. Such an extension $V=\oplus_{g\in G}V^g$
graded by a finite abelian group $G$
can be characterized the properties listed below. We first obtain
the complete reducibility of a $G$-graded module $M=\oplus_{g\in G}M^g$ with
$M^0$ being an irreducible $V^0$-module. Then we study the complete
reducibility
of a canonical class of $V$-modules.
We apply our results to some special cases.

{}From the last section (Corollaries 3.13, 3.14 and 3.15) under certain
conditions we obtain
vertex operator (super)algebras satisfy the following conditions:

(1) $V=\oplus_{g\in G} V^{g}$, where $G$ is a finite abelian group.

(2) $V^{0}$ is a simple rational vertex operator subalgebra and each $V^{g}$ is
a
simple current $V^{0}$-module.

(3) $a_{n}V^{h}\subseteq V^{g+h}$ for any $a\in V^{g}, n\in {\Z}, g,h\in G$.

(4) For any subgroup $H$ of $G$, $V^{(H)}=: \oplus_{h\in H}V^{(h)}$
is a simple vertex (operator) algebra and
$V^{(g+H)}=:\oplus_{h\in H}V^{(h+g)}$ for any $g\in G$ is a simple current
$V^{(H)}$-module.

(5) Let $M^{0}$ be any irreducible $V^{0}$-module which is a $V^{0}$-submodule
of some $V$-module $W$. Then there is a $V$-module $M=\oplus_{g\in G} M^{g}$
satisfying
the condition: $a_{n}M^{h}\subseteq M^{g+h}$ for any $a\in V^{g}, n\in {\Z},
g,h\in G$.

\begin{prop}\label{p4.1} {\it Let $M=\oplus_{g\in G}M^{g}$ be a $V$-module
satisfying
the following conditions: (a) $M^{0}$  is an irreducible $V^{0}$-module.
(b) $a_{n}M^{h}\subseteq M^{g+h}$ for any $a\in V^{g}, n\in {\Z}, g,h\in G$.\\
Then $M$ is a direct sum of irreducible $V$-modules.}
\end{prop}

First we prove the following two special cases:

\begin{lem}\label{l4.2} {\it Let $M=\oplus_{g\in G}M^{g}$ be a $V$-module
satisfying
the conditions (a) and (b) of Proposition 4.1. Suppose $M^{g}$ and $M^{h}$ are
not isomorphic $V^{0}$-modules for $g\ne h$. Then $M$ is an irreducible
$V$-module.}
\end{lem}

{\bf Proof.} Let $M_{1}$ be any nonzero $V$-submodule of $M$. We must show that
$M=M_{1}$.
Since each $M^{g}$ generates $M$ by the action of $V$, it suffices to prove
that $M_{1}$
contains some $M^{g}$. Because $M$ is a direct sum of irreducible
$V^{0}$-modules, $M_{1}$
is also a direct sum of irreducible $V^{0}$-modules.
For any $g\in G$, let $P_{g}$ be the projection of $M$ onto $M^{g}$. Then
$P_{g}$ is a $V^{0}$-homomorphism. Let $W$ be any irreducible $V^{0}$-submodule
of $M_{1}$.
Then the restriction of $P_{g}$ to $W$ is either zero or a $V^{0}$-isomorphism
onto $M^{g}$. Since
$M^{g}$ and $M^{h}$ are not $V^{0}$-isomorphic for $g\ne h$, there is $g\in G$
such that
$P_{g}(W)=M^{g}$ and $P_{h}(W)=0$ for
$h\ne g$. Therefore $W=M^{g}$ for some $g\in G$. Thus $M_{1}$ contains some
$M^{g},$ as required. $\;\;\;\;\Box$

\begin{lem}\label{l4.3} {\it Let $M=\oplus_{g\in G}M^{g}$ be a $V$-module
satisfying conditions (a) and (b) of Proposition 4.1. Suppose that $G$ is a
cyclic group and that
all $M^{g}$ ($g\in G$) are isomorphic irreducible $V^{0}$-modules.
Then $M$ is a direct sum of $|G|$ irreducible $V$-modules, each of which is
isomorphic
 to $M^{0}$ as a $V^{0}$-module.}
\end{lem}

{\bf Proof.} For any $g\in G$, let $f_{g}$ be a fixed $V^{0}$-isomorphism from
$M^{0}$
onto $M^{g}$. If $g=0$, we choose $f_{0}=Id_{M^{0}}$.
We shall extend $f_{g}$ to be a $V$-automorphism  of $M$.
Let $h\in G$. Then $(M^{h}, Y(\cdot,z))$ is a tensor product
for $(V^{h},M^{0})$ by Corollary 2.9.
Since $Y(\cdot,z)f_{g}$ is an intertwining operator of type
$\left(\begin{array}{c}M^{g+h}\\V^{h},M^{0}\end{array}\right)$, there is a
(unique)
$V^{0}$-homomorphism $f_{g,h}$ from $M^{h}$ to $M^{g+h}$ satisfying
\begin{eqnarray}
Y(a,z)f_{g}(u)=f_{g,h}Y(a,z)u\;\;\;\mbox{for any }a\in V^{h}, u\in M^{0}
\end{eqnarray}
(see Definition 2.3 for the tensor product of modules).
If $h=0$, we have $f_{g,0}=f_{g}$.
Then we extend $f_{g}$ to be a  linear endomorphism of $M$ by defining:
\begin{eqnarray}
f_{g}(u)=f_{g,h}(u)\;\;\;\mbox{for any }h\in G, u\in M^{h}\subseteq M.
\end{eqnarray}
Thus
\begin{eqnarray}
Y(a,z)f_{g}(u)=f_{g}Y(a,z)u\;\;\;\mbox{for any }a\in V, u\in M^{0}.
\end{eqnarray}
Let $a,b\in V,g\in G, u\in M^{0}$ and let $k$ be a positive integer such that
the
following associativities hold:
\begin{eqnarray}
&
&(z_{0}+z_{2})^{k}Y(a,z_{0}+z_{2})Y(b,z_{2})u=(z_{0}+z_{2})^{k}Y(Y(a,z_{0})b,z_{2})u,\\
& &(z_{0}+z_{2})^{k}Y(a,z_{0}+z_{2})Y(b,z_{2})f_{g}u
=(z_{0}+z_{2})^{k}Y(Y(a,z_{0})b,z_{2})f_{g}u.
\end{eqnarray}
Then by (4.3)-(4.5) we have:
\begin{eqnarray}
& &(z_{0}+z_{2})^{k}f_{g}Y(a,z_{0}+z_{2})Y(b,z_{2})u\nonumber\\
&=&(z_{0}+z_{2})^{k}f_{g}Y(Y(a,z_{0})b,z_{2})u\nonumber\\
&=&(z_{0}+z_{2})^{k}Y(Y(a,z_{0})b,z_{2})f_{g}(u)\nonumber\\
&=&(z_{0}+z_{2})^{k}Y(a,z_{0}+z_{2})Y(b,z_{2})f_{g}(u)\nonumber\\
&=&(z_{0}+z_{2})^{k}Y(a,z_{0}+z_{2})f_{g}(Y(b,z_{2})u).
\end{eqnarray}
Multiplying by $(z_{0}+z_{2})^{-k}$ we obtain
\begin{eqnarray}
f_{g}Y(a,z_{0}+z_{2})Y(b,z_{2})u=Y(a,z_{0}+z_{2})f_{g}(Y(b,z_{2})u).
\end{eqnarray}
Thus
\begin{eqnarray}
f_{g}Y(a,z_{1})Y(b,z_{2})u=Y(a,z_{1})f_{g}(Y(b,z_{2})u).
\end{eqnarray}
Since $V\cdot M^{0}=M$, we get
\begin{eqnarray}
f_{g}Y(a,z_{1})v=Y(a,z_{1})f_{g}v\;\;\;\;\mbox{for any }a\in V, v\in M.
\end{eqnarray}
That is, $f_{g}$ is a $V$-endomorphism of $M$.

For any $g,h\in G$, both $f_{g+h}$
and $f_{g}f_{h}$ are $V^{0}$-homomorphisms from $M^{0}$ to $M^{g+h}$. Since
$M^{0}$ is
an irreducible $V^{0}$-module, $f_{g+h}$ is a constant multiple of
$f_{g}f_{h}$.
That is, there is $A(g,h)\in {\C}^{*}$ such that $f_{g+h}=A(g,h)f_{g}f_{h}$
{}from $M^{0}$ to $M^{g+h}$. Since each $f_{g}$ is a $V$-endomorphism of $M$
and
$M^{0}$ generates $M$ by $V$,
$f_{g+h}=A(g,h)f_{g}f_{h}$ holds on $M$. It is clear that $A(g,h)$ is a
$2$-cocycle.

Since $G$ is cyclic, let $G=\langle g\rangle$ with $o(g)=k$. Since
$f_{g}^{k}$ is a $V^{0}$-endomorphism of $M^{0}$ and $M^{0}$ is an irreducible
$V^{0}$-module, there is a complex number $\alpha$ such that
$f_{g}^{k}(u)=\alpha u$ for
any $u\in M^{0}$. Then we modify $f_{g}$ by multiplying a $k$-th root of
$\alpha$, we
have: $f_{g}^{k}(u)=u$ for any $u\in M^{0}$. Since $f_{g}$ commutes with all
vertex
operators $Y(a,z)$ for $a\in V$ and $M^{0}$ generates $M$ by $V$, we have
$f_{g}^{k}=Id_{M}$.
Using $f_{g}$ we obtain a representation of $G$ on $M$. For any nonzero $u\in
M^{(0)}$,
${\C}u\oplus {\C}f_{g}u\oplus\cdots\oplus \C f_{g}^{k-1}u$ is isomorphic to the
regular
representation of $G$. For any character $\chi\in \hat{G}$, let $M(\chi)$ be
the $\chi$-homogeneous subspace of $M.$ Then
$M=\oplus _{\chi \in \hat{G}}M(\chi)$ and $M(\chi)\ne 0$ for any $\chi\in
\hat{G}$.
Since $G$ commutes with all vertex
operators $Y(a,z)$ for $a\in V$, each
$M(\chi)$ is a $V$-module. Since $M(\chi)\ne 0$ for any $\chi\in \hat{G}$ and
$M$ is
a direct sum of $|G|$ irreducible $V^{0}$-modules, then each $M(\chi)$ must be
an
irreducible $V$-module, so that it is also an irreducible
$V^{0}$-module.$\;\;\;\;\Box$

{\bf Proof of Proposition 4.1.} We are going to prove Proposition 4.1 by using
induction on $|G|$.
If $|G|=1$, there is nothing to prove. Suppose that Proposition 4.1 is true for
any finite abelian
group with less than $n$ elements. Suppose that $|G|=n$ with $n>1$.
Let $H$ be a subgroup of $G$ of prime order. Set $V^{(H)}=\oplus_{h\in
H}V^{(h)}$ and
$V^{(g+H)}=\oplus_{h\in H}V^{(g+h)}$ for $g\in G$. Then $V^{(H)}$ is a simple
vertex operator
algebra and each $V^{(g+H)}$ is a simple current for $V^{(H)}$.
Similarly set $M^{(H)}=\oplus _{h\in H}M^{(h)}$. Let $H_{0}$ be the subset of
$H$ consisting of
$h$ such that $M^{(h)}$ is isomorphic to $M^{(0)}$. Then it is clear that
$H_{0}$ is a subgroup
of $H$. Consequently, either $H_{0}=H$ or $H_{0}=0$. By Lemmas 4.2 and 4.3,
$M^{(H)}$ is a direct
sum of irreducible $V^{(H)}$-modules. Let $M^{(H)}=W_{1}\oplus \cdots\oplus
W_{m}$, where $W_{j}$ are irreducible $V^{(H)}$-modules. Let $M^j$ be the
$V$-module
generated by $W_j.$ Denote the span of $\{a_nw|a\in V^g,n\in\Z,w\in W_j\}$
by $V^g W_j$ for $g\in G.$   Then from the proof of Lemmas 4.2 and 4.3,
$$M_j=\sum_{g\in G}V^g W_j=\oplus_{k\in G/H}V^kW_j$$
and $M=\sum_{j=1}M^j.$ Then $V^{(H)}, G/H, M_{j}$
satisfy the assumptions of Proposition 4.1. By the inductive assumption, each
$M_{j}$ is a direct sum of
irreducible $V$-modules, hence so too is $M$.$\;\;\;\;\Box$

\bt{t4.5} {\it Suppose that $V^{0}$ is rational and that for any irreducible
$V^{0}$-module $W^{0}$, if there is a
$V$-module $W$ such that $W^{0}$ is a $V^{0}$-submodule of $W$, then $W^{0}$
can be lifted
to be a $V$-module $M=\oplus_{g\in G}M^{g}$ with $M^{0}=W^{0}$. Then $V$ is
rational.}
\et

{\bf Proof.} Let $W$ be any $V$-module. Then $W$ is a completely reducible
$V^{0}$-module.
Thus it suffices to prove that any irreducible $V^{0}$-submodule $W^{0}$ of $W$
generates
 a completely reducible $V$-submodule of $W$. By assumption, there is a
$V$-module
$M=\oplus_{g\in G}M^{g}$ such that
\begin{eqnarray}
M^{0}=W^{0} ,\;\;a_{n}M^{h}\subseteq M^{g+h}\;\;\;\mbox{for any }a\in V^{g},
 n\in {\Z},g,h\in G.
\end{eqnarray}
{}From Proposition 4.1, $M$ is a completely reducible $V$-module. So it is
sufficient to
prove that the $V$-submodule $\langle W^{0}\rangle$ of $W$ generated by $W^{0}$
is a $V$-homomorphism image of $M$. It is easy to see that $(M, Y(\cdot,z))$ is
a tensor product
of $V^{0}$-modules $V$ with $M^{0}=W^{0}$. Therefore there is a
$V^{0}$-homomorphism $f$
{}from $M$ to $W$ such that
\begin{eqnarray}
Y_{W}(a,z)u=fY_{M}(a,z)u\;\;\;\mbox{for any }a\in V, u\in M^{0}=W^{0}.
\end{eqnarray}
Using the same argument used in the proof of Lemma 4.3 we obtain:
\begin{eqnarray}
Y_{W}(a,z)f(u)=fY_{M}(a,z)u\;\;\;\mbox{for any }a\in V, u\in M.
\end{eqnarray}
Then $f$ is a $V$-homomorphism from $M$ to $W$. Then $\langle W^{0}\rangle$ is
a
completely reducible $V$-module. Therefore $M$ is a completely reducible
$V$-module.
$\;\;\;\;\Box$.

\bt{t4.6} {\it Suppose that $V^{0}$ is rational and that for any $g\in G$,
there is an
$h_{g}\in V_{1}$
satisfying condition (2.18) and such that $V^{g}$ is $V^{0}$-isomorphic to
$(V^{0}, Y(\Delta(h_{g},z)\cdot,z))$.
Then $V$ is rational.}
\et

{\bf Proof.} This follows from Theorem \ref{t4.5} and Corollaries 3.14-15
immediately.\hspace{1cm}$\Box$

\section{Applications to affine Lie algebras}

This section is devoted to the study of affine Kac-Moody algebras and their
representations. It is well known that $L(l\Lambda_0)$ is a vertex
operator algebra (see the definition below) and any weak
module which is truncated below is a direct sum of standard modules of level
$l$ (cf. [DL] and [FZ]).
In this section we improve this result by showing
that under a mild assumption, any weak module is a direct sum of standard
modules of level $l.$ Then we discuss the simple currents for vertex operator
algebras $L(l\Lambda_0)$ and various extensions of $L(l\Lambda_0)$ as
applications
of results obtained in previous sections.

Let ${\bf g}$ be a finite-dimensional simple Lie
algebra with a fixed Cartan subalgebra $H$ and let $\{e_{i}, f_{i},
\alpha^{\vee}_{i}|i=1,\cdots, n\}$ be the Chevalley generators. Let
$(\cdot,\cdot)$ be the
normalized Killing form on ${\bf g}$ such that the square norm of the longest
root is $2$.
Let
${\Q}={\Z}\alpha_{1}\oplus\cdots\oplus {\Z}\alpha_{n}$, where
$\alpha_{1},\cdots,\alpha_{n}$ are all simple roots. Notice that
$\alpha^{\vee}_{1},\cdots,\alpha^{\vee}_{n}$ form a basis for $H$.  Let
$h_{i}\in H$ such that
$\alpha_{i}(h_{j})=\delta_{i,j}$ for $i,j=1,\cdots,n$.
Let
$\displaystyle{\theta=\sum_{i=1}^{n}a_{i}\alpha_{i}}$ be the highest
positive root.  Let $\lambda_{i}$ ($i=1,\cdots,n$) be the fundamental
weights for ${\bf g}$ (cf. [H]).  A dominant integral weight $\lambda$
is called a {\em minimal} weight (cf. [H]) if there is no dominant
integral weight $\gamma$ satisfying $\lambda-\gamma\in {\Q}_{+}$.
Then $\lambda_{i}$ is minimal if and only if $a_{i}=1,$ and all minimal
dominant integral weights are given as follows (cf. [H]):
\begin{eqnarray}
 A_{n}:& &\lambda_{1},\cdots, \lambda_{n}\nonumber\\
B_{n}: & &\lambda_{n}\nonumber\\
C_{n}:& & \lambda_{1}\nonumber\\
D_{n}: & & \lambda_{1},\lambda_{n-1},\lambda_{n}\nonumber\\
E_{6}:& &\lambda_{1},\lambda_{6}\nonumber\\
E_{7}: & &\lambda_{7}.
\end{eqnarray}

Let $\tilde{{\bf g}}$ be the affine Lie algebra [K] with Chevalley generators
$\{e_{i}, f_{i}, \alpha^{\vee}_{i}|i=0,\cdots, n\}$. Then each
$\lambda_{i}$ for $1\le i\le n$ is naturally extended to a fundamental weight
$\Lambda_{i}$
for $\tilde{{\bf g}}$.
Let $\Lambda_{0}$ be the fundamental weight for $\tilde{{\bf g}}$
defined by $\Lambda_{0}(\alpha^{\vee}_{i})=\delta_{i,0}$ for $0\le i\le n$ (cf.
[K]). Then
$\Lambda_{i}$ is of level one if and only if $a^{\vee}_{i}=1$ (see
[K] for the definition of $a^{\vee}_i$).
Let $\lambda\in H^{*}$ and let $\ell$ be any complex number. Then we denote by
$L(\ell,\lambda)$ the highest weight $\tilde{{\bf g}}$-module of level
$\ell$ with
lowest weight $\lambda.$
 It is well known (cf. [DL], [FZ], [Li1]) that $L(\ell,0)$
is a vertex operator algebra. One can identify ${\bf g}$ as a subspace of
$L(\ell,0)$
 through the linear map $\phi: u\mapsto u_{-1}{\bf 1}$.
Using this formulation, it was proved in [Li4] that if $\lambda_{i}$ is
minimal, then
for any $\ell$, $L(\ell\Lambda_{i})$ (or $L(\ell,\lambda_{i})$) is an
irreducible
(weak) $L(\ell,0)$-module and it is a simple current.
The following proposition is proved in [Li4]. Since the information obtained
is useful (see Remarks 5.2 and 5.3), we repeat the short proof here.

\begin{prop}\label{p5.1} {\it Suppose that $\lambda_{i}$ is minimal.
Let $\ell$ be  any complex number which is not equal to $-\Omega$,
where $\Omega$ is the dual Coxeter number of ${\bf g}$. Then
 $L(\ell\Lambda_{i})$ is isomorphic to $(L(\ell,0), Y(\Delta(h_{i},z)\cdot,z))$
 as a $\tilde{{\bf g}}$-module.
Consequently, $L(\ell\Lambda_{i})$ is a simple current.}
\end{prop}

{\bf Proof.} Note that $\theta (h_i)=a_{i}=1$. By definition we have:
\begin{eqnarray}
& &\Delta(h_i,z)\alpha^{\vee}_{j}=\alpha^{\vee}_{j}+\ell \delta_{i,j}z^{-1},\;
\Delta(h_i,z)e_{i}=ze_{i},\;
\Delta(h_i,z)f_{i}=z^{-1}f_{i},\\
& &\Delta(h_i,z)e_{j}=e_{j},\;
\Delta(h_i,z)f_{j}=f_{j},\;\Delta(h_i,z)f_{\theta}=z^{-1}f_{\theta}
\;\;\mbox{ for }j\ne i.
\end{eqnarray}
In other words, the corresponding automorphism $\psi$ of
$U(\tilde{{\bf g}})$ or $U(L(\ell,0))$ satisfies the following conditions:
\begin{eqnarray}
& &\psi(\alpha^{\vee}_{i}(n))=\alpha^{\vee}_{i}(n)+\delta_{n,0}\ell,\;
\psi(e_{i}(n))=e_{i}(n+1),\;
\psi(f_{i}(n))=f_{i}(n-1);\\
& &\psi(\alpha^{\vee}_{j}(n))=\alpha^{\vee}_{j}(n),\;\psi(e_{j}(n))=e_{j}(n),\;
\psi(f_{j}(n))=f_{j}(n)
\;\;\;\mbox{for }j\ne i, n\in {\Z},
\end{eqnarray}
and
\begin{eqnarray}
\psi (f_{\theta}(n))=f_{\theta}(n-1)\;\;\;\mbox{for }n\in {\Z}.
\end{eqnarray}
Then the vacuum vector ${\bf 1}$ in $(L(\ell,0), Y(\Delta(h_i,z)\cdot,z))$ is a
highest weight vector of weight $\ell\Lambda_{i}$. Thus $(V,
Y(\Delta(h_i,z)\cdot,z))$ is isomorphic to $L(\ell\Lambda_{i})$ as a
$\tilde{{\bf g}}$-module. By Proposition 2.12, $L(\ell\Lambda_{i})$ is a
simple current. $\;\;\;\;\Box$

\begin{rem}\label{r5.2}
It has been shown  [FG] by calculating the four point functions that if
$\lambda_{i}$ is a minimal weight and $\ell$ is a positive
integer, then $L(\ell\Lambda_{i})$ is a simple current. Moreover
if ${\bf g}$ is of type $E_{8}$, $L(\Lambda_{7})$ is a simple current of level
$2$ which is not isomorphic to $L(2,0).$
 It has also been proved [F] that these are
the all simple currents. In string theory, simple currents are useful for
constructing
modular invariants.
\end{rem}

\begin{rem}\label{r5.3} From the proof of Proposition 5.1 we see that
the vacuum ${\bf 1}$ becomes a lowest weight vector of $L(\ell\Lambda_{i})$.
The lowest weight of $L(\ell\Lambda_{i})$ is ${\ell\over 2}\langle
h_i,h_i\rangle$
because $L(0)$ acts on
$L(\ell\Lambda_{i})$ $(=L(\ell\Lambda_{0}))$ as
$L(0)+h_{i}(0)+{\ell\over 2}\langle h_{i},h_{i}\rangle$.
In general, if $h\in H$ satisfies $\alpha_{i}(h)\in {\Z}$ for $i=1,\cdots,n$,
it follows from the proof of Proposition 5.1 that the vacuum vector is a
lowest
weight vector for $\tilde{{\bf g}}$ if and only if either $h=0$ or $h$
corresponds to a
minimal weight. Then for some $h\in H$,
$(L(\ell,0), Y(\Delta(h,z)\cdot,z))$ might not be a highest weight $\tilde{{\bf
g}}$-module.
\end{rem}

Next we shall prove
that if $\ell$ is a positive integer, then for any $h\in H$ satisfying
$\alpha_{i}(h)\in {\Z}$ for $i=1,\cdots,n$ and
for any $L(\ell,0)$-module $(M,Y_{M}(\cdot,z))$,
 $(M,Y_{M}(\Delta(h,z)\cdot,z))$ is an ordinary
$L(\ell,0)$-module. {\em From now on, we assume that ${\bf g}$ is a fixed
finite-dimensional
simple Lie algebra and $\ell$ is a fixed positive integer.}

The following lemma easily follows from Proposition 13.16 in [DL] (see [Li1] or
[MP] for a proof).

\begin{lem}\label{l5.4} {\it Let $M$ be any weak $L(\ell,0)$-module and let
$e\in {\bf g}_{\alpha}$, where $\alpha$ is any root of ${\bf g}$. Then
 $Y_{M}(e,z)^{\ell+1}=0$ if $\alpha$ is a long root and
$Y_{M}(e,z)^{3\ell+1}=0$ for
any $\alpha$.}
\end{lem}

\begin{lem}\label{l5.5} {\it Let $M$ be any nonzero weak
$L(\ell,0)$-module on which $t{\C}[t]\otimes H$ acts locally nilpotently.
Then  $M$ contains a standard $\tilde{{\bf g}}$-module of level
$\ell$.}
\end{lem}

{\bf Proof.} Set $\tilde{{\bf g}}_{+}=t{\C}[t]\otimes {\bf g}$.
We define $M^{0}=\{u\in M| \tilde{{\bf g}}_{+}u=0\}$. Since
$[{\bf g},\tilde{{\bf g}}_{+}]\subseteq \tilde{{\bf g}}_{+}$,
${\bf g}M^{0}\subseteq M^{0}$. Let $0\ne e\in {\bf g}_{\theta}$.
Applying $Y_{M}(e,z)^{\ell+1}$ to
$M^{0}$ and extracting the coefficient of $z^{-\ell-1}$, we obtain
$e(0)^{\ell+1}M^{0}=0$.
By Proposition 5.1.2 in [Li1] (see also [KW]), $M^{0}$ is a direct sum of
finite-dimensional
irreducible ${\bf g}$-modules. If $M^{0}\ne 0$, let $u$ be a highest weight
vector
for ${\bf g}$ in $M^{0}$. Then  $u$ is a highest weight vector for $\tilde{{\bf
g}}$.
Extracting the constant from $Y_{M}(e,z)^{\ell+1}u=0$
we obtain $e(-1)^{\ell+1}u=0$.
Then $u$ generates a standard $\tilde{{\bf g}}$-module. So it suffices to prove
that
$M^{0}$ is nonzero.
For any $u\in M$, it follows from the definition of a weak $L(\ell,0)$-module
that $\tilde{{\bf g}}_{+}u$ is finite-dimensional. Let $u$ be a nonzero vector
 of $M$ such that $\tilde{{\bf g}}_{+}u$ has minimal dimension. If the
minimal dimension is zero, then we are done. Suppose the minimal dimension is
not zero. Let $k$  be an integer such that ${\bf g}(n)u=0$ for $n>k$ and that
there is a nonzero $a\in {\bf g}$ such that $a(k)u\ne 0$. Without loss
generality
we may assume that $a\in {\bf g}_{\alpha}$ for some root $\alpha$ or $a\in H$.
By assumption,
$k>0$.  If $a\in {\bf g}_{\alpha}$ for some root $\alpha$, then
$Y_{M}(a,z)^{3\ell+1}u=0$ implies
that $a(k)^{3\ell+1}u=0$. If $a\in H$, by assumption $a(k)$ locally nilpotently
acts on $u$.
Therefore, there is a
nonnegative integer $r$ such that $a(k)^{r}u\ne 0$ and $a(k)^{r+1}u=0$. Let
$v=a(k)^{r}u$. If $b(n)u=0$ for some $b\in {\bf g}, n>1$, then $b(n)v=0$. Then
$\dim \tilde{{\bf g}}_{+}v\le \dim\tilde{{\bf g}}_{+}u-1,$
contradiction. The proof is complete.$\;\;\;\;\Box$

\begin{prop}\label{p5.6} {\it Let $M$ be a
weak $L(\ell\Lambda_{0})$-module on which $t{\C}[t]\otimes H$ acts locally
nilpotently. Then
$M$ is a direct sum of standard $\tilde{{\bf g}}$-modules of level $\ell$.}
\end{prop}

{\bf Proof.} Let $M_{1}$ be a direct sum of standard
submodules of $M$. We have to prove that $M=M_{1}$. Otherwise,
the quotient module $\bar{M}=M/M_{1}$ is not zero. By lemma 5.5, there is a
standard submodule of $\bar{M}$, say $W/M_{1}$, where $W$ is a submodule of
$M$ containing $M_{1}$. It follows from Theorem 10.7 in [K] that $W$ is a
direct sum of
standard modules, so that $W=M_{1},$ contradiction.$\;\;\;\Box$

\bc{c5.7} {\it Let $h\in H$
such that $\alpha (h)\in {\Z}$ for any root $\alpha$ of ${\bf g}$ and let
$(M,Y_{M}(\cdot,z))$  be any $L(\ell,0)$-module. Then
$(M, Y(\Delta(h,z)\cdot,z))$ is also an  $L(\ell,0)$-module.}
\ec

{\bf Proof.} Since any $L(\ell,0)$-module $M$ is a direct sum of finitely many
standard
$\tilde{{\bf g}}$-modules and a direct sum of finitely many $L(\ell,0)$-modules
is a module, it
suffices to prove the corollary for an irreducible module $M$.
Since $\Delta(h,z)$ is invertible, it is clear that
$(M, Y(\Delta(h,z)\cdot,z))$  is still irreducible as a $\tilde{{\bf
g}}$-module.
It follows from the proof of Proposition 2.15 that
$(M, Y(\Delta(h,z)\cdot,z))$ is still a direct sum of highest weight modules
for
the Heisenberg Lie algebra $\hat{H}$. By Proposition 5.6, $(M,
Y(\Delta(h,z)\cdot,z))$
is a direct sum of standard $\tilde{{\bf g}}$-modules of level $\ell$.
Consequently, $M$ is a
standard $\tilde{{\bf g}}$-module of level $\ell$. The proof is
complete.$\;\;\;\;\Box$

\begin{rem}\label{r5.8} If we just consider the action of the affine Lie
algebra
$\tilde{{\bf g}}$,
it is easy to see that the transition from $L(\ell,0)$ to $L(\ell,\lambda_{i})$
is due to an
Dynkin diagram automorphism of $\tilde{{\bf g}}$.
The automorphism group of the Dynkin diagram of
the affine Lie algebra is commonly called the {\em outer automorphism group}.
To a certain
extent we have
realized them explicitly as  ``inner automorphisms'' in terms of exponentials
of certain
elements of $\tilde{{\bf g}}$.
\end{rem}

\begin{rem}\label{r5.9}
It is very special for $E_{8}$ that there is a simple
current other than the vacuum representation when $\ell =2$, but there is no
outer automorphism. Let $h\in H$ be the element uniquely determined by
$\alpha_{i}(h)=\delta_{i,7}$ for $1\le i\le 7$. By Corollary 5.7,
$(L(2,0),Y(\Delta(h,z)\cdot,z))$ is still a
standard module and it is a simple current. It is interesting to ask what this
module is.
If this module is $L(2,\lambda_{7})$, then all simple currents can be
constructed in terms of
$\Delta(h,z)$.
\end{rem}

Let $L$ be the ${\Z}$-span of all minimal weights of ${\bf g}$. Then it follows
{}from
Theorem 3.21 that
for any positive integral level $\ell$, the direct sum of all simple currents
of
$L(\ell,0)$ is an abelian intertwining algebra.

 \bt{t5.10} {\it Let ${{\bf g}}$ be a finite-dimensional simple Lie
algebra, but not of type $E_{8}$ with a fixed Cartan subalgebra $H$ and let
$\ell$ be any positive integer. Then the direct sum of all (simple currents)
$L(\ell,0)$-modules $L(\ell,\lambda_{i})$ for all minimal weights $\lambda_{i}$
is an
abelian intertwining algebra.}
\et

{\bf Example 5.11} Let ${\bf g}$ be of type $A_{n}$. From  [H] we have:
\begin{eqnarray}
\lambda_{i}&=&\frac{1}{ n+1}\left((n-i+1)\alpha_{1}+2(n-i+1)\alpha_{2}+\cdots+
(i-1)(n-i+1)\alpha_{i-1}\right)\nonumber\\
& &+\frac{1}{
n+1}\left(i(n-i+1)\alpha_{i}+i(n-i)\alpha_{i+1}+\cdots+i\alpha_{n}\right).
\end{eqnarray}
Then
$$h_{i}=\frac{1}{
n+1}\left((n-i+1)\alpha_{1}^{\vee}+2(n-i+1)\alpha_{2}^{\vee}+\cdots+
(i-1)(n-i+1)\alpha_{i-1}^{\vee}\right)$$
$$\ \ \ \ \ +\frac{1}{
n+1}\left(i(n-i+1)\alpha_{i}^{\vee}+i(n-i)\alpha_{i+1}^{\vee}+\cdots+
i\alpha_{n}^{\vee}\right).$$
By a simple calculation we get
$(h_{i},h_{i})=\frac{i(n+1-i)}{n+1}$ for $1\le i\le n$.
Let $L=L_{1}={\Z}h_{1}$.
Notice that (from (\ref{3.7}))
\begin{eqnarray}
\langle h,h'\rangle {\bf 1}=h(1)h'=h(1)h'(-1){\bf 1}=\ell (h,h')
\end{eqnarray}
for $h,h'\in H$. Then
$\langle h_{i},h_{i}\rangle=\frac{\ell i(n+1-i)}{n+1}$. By definition, $P$ is
the root lattice
with the bilinear form $\ell (\cdot,\cdot)$.
If $\ell\in 2(n+1){\Z}_{+}$,
$L_{1}$ satisfies condition (3.56).  Suppose that $L_{0}=kL_{1}$ for some
positive integer $k$.
By Corollaries 3.13 and 2.9,
$L(l,0)\oplus \oplus_{i=0}^{k-1}L(\ell,\lambda_{1})^{{\boxtimes}i}$ is a vertex
operator algebra if
$\ell\in 2(n+1){\Z}_{+}$. In fact, $L_{0}=(n+1)L_{1}$. By Theorem \ref{t4.6},
it is rational.

For other types, there are similar arguments, so we just briefly state the
result.

{\bf Example 5.12.} If ${\bf g}$ is of type $B_{n}$, then (cf. [H])
\begin{eqnarray}
\lambda_{n}={1\over 2}(\alpha_{1}+2\alpha_{2}+\cdots+n\alpha_{n}).
\end{eqnarray}
Thus $h_{n}={1\over
2}(\alpha^{\vee}_{1}+2\alpha^{\vee}_{2}+\cdots+n\alpha^{\vee}_{n})$.
Let $L=L_{1}={\Z}h_{n}$. Since
$L(\ell,\lambda_{n})$ is isomorphic to its contragredient module, we have a
nonzero intertwining
operator of type
$\left(\begin{array}{c}L(\ell,0)\\L(\ell,\lambda_{n})
L(\ell,\lambda_{n})\end{array}\right)$.
Then $2h_{n}\in L_{0}$. Thus $L_{0}=2L_{1}$.
Since $(h_{n},h_{n})=1$, we have $\langle \lambda_{n},\lambda_{n}\rangle
=\ell$.
Then by Corollary 3.13 and Theorem \ref{t4.6}, $L(\ell,0)\oplus
L(\ell,\lambda_{n})$
is a rational vertex operator algebra for $\ell$ is even and
$L(\ell,0)\oplus L(\ell,\lambda_{n})$
is a rational vertex operator superalgebra if $\ell$ is odd. If $\ell=1$, the
lowest
weight of $L(1,\lambda_{n})$
is ${1\over 2}$. It follows from the super Jacobi identity that
$L(1,\lambda_{n})_{{1\over 2}}$  generates a Clifford algebra.
This result explains why one can realize $\tilde{B}_{n}$ [LP] in terms of the
representations
of a Clifford algebra.

{\bf Example 5.13.} If ${\bf g}$ is of type $C_{n}$, then
\begin{eqnarray}
\lambda_{1}=\alpha_{1}+\cdots+\alpha_{n-1}+{1\over 2}\alpha_{n}.
\end{eqnarray}
Thus $h_{1}=\alpha^{\vee}_{1}+\cdots+\alpha^{\vee}_{n-1}+{1\over
2}\alpha^{\vee}_{n}$.
Since $(h_{1},h_{1})=1$, we have: $\langle h_{1},h_{1}\rangle=\ell$.
Thus $L(\ell,0)\oplus L(\ell,\lambda_{n})$ if $\ell$ is even
is a rational vertex operator algebra and
$L(\ell,0)\oplus L(\ell,\lambda_{n})$
is a rational vertex operator superalgebra if $\ell$ is odd.

{\bf Example 5.14.} If ${\bf g}$ is of type $D_{n}$, then
\begin{eqnarray}
\lambda_{1}&=&\alpha_{1}+\cdots+\alpha_{n-2}+{1\over
2}(\alpha_{n-1}+\alpha_{n}),\nonumber\\
\lambda_{n-1}&=&{1\over 2}\left(\alpha_{1}+2\alpha_{2}+\cdots+(n-2)\alpha_{n-2}
+{n\over 2}\alpha_{n-1}+\frac{n-2}{2}\alpha_{n}\right),\nonumber\\
\lambda_{n}&=&{1\over 2}\left(\alpha_{1}+2\alpha_{2}+\cdots+(n-2)\alpha_{n-2}
+\frac{n-2}{2}\alpha_{n-1}+\frac{n}{2}\alpha_{n}\right).
\end{eqnarray}
Then
\begin{eqnarray}
h_{1}&=&\alpha^{\vee}_{1}+\cdots+\alpha^{\vee}_{n-2}+{1\over
2}(\alpha^{\vee}_{n-1}
+\alpha^{\vee}_{n}),\nonumber\\
h_{n-1}&=&{1\over
2}\left(\alpha^{\vee}_{1}+2\alpha^{\vee}_{2}+\cdots+(n-2)\alpha^{\vee}_{n-2}
+{n\over
2}\alpha^{\vee}_{n-1}+\frac{n-2}{2}\alpha^{\vee}_{n}\right),\nonumber\\
h_{n}&=&{1\over
2}\left(\alpha^{\vee}_{1}+2\alpha^{\vee}_{2}+\cdots+(n-2)\alpha^{\vee}_{n-2}
+\frac{n-2}{2}\alpha^{\vee}_{n-1}+\frac{n}{2}\alpha^{\vee}_{n}\right).
\end{eqnarray}
By a simple calculation we obtain:
\begin{eqnarray}
(h_{1},h_{1})=1,\; (h_{n-1},h_{n-1})=(h_{n},h_{n})={n\over 4}.
\end{eqnarray}
Then
\begin{eqnarray}
\langle h_{1},h_{1}\rangle=\ell,\; \langle h_{n-1},h_{n-1}\rangle
=\langle h_{n},h_{n}\rangle={n\ell\over 4}.
\end{eqnarray}
Let $L=L_{1}={\Z}h_{1}\oplus {\Z}h_{n-1}\oplus {\Z}h_{n}$. It is not difficult
to see that
$L_{0}=2L_{1}$.
Thus
$L(\ell,0)\oplus L(\ell,\lambda_{1})\oplus L(\ell,\lambda_{n-1})\oplus
L(\ell,\lambda_{n})$
is an abelian intertwining algebra with ${\Z}_{2}\times {\Z}_{2}$ as its
grading group.
(This has been proved in [DL].)
Furthermore, $L(\ell,0)\oplus L(\ell,\lambda_{1})$ is a rational vertex
operator superalgebra.
Another special case is when $n\ell$ is divisible by
$8$ [DM]: $L(\ell,0)\oplus L(\ell,\lambda_{n-1})$ and $L(\ell,0)\oplus
L(\ell,\lambda_{n})$
are (holomorphic) vertex operator algebras.

\end{document}